\documentclass{amsart}
\usepackage{preamble_w}

\newcommand{\setOfReals}{\mathbb{R}}
\newcommand{\setOfNaturals}{\mathbb{N}}
\newcommand{\setOfNonnegativeIntegers}{{\mathbb{N}_0}}
\newcommand{\setOfPositiveReals}{{\setOfReals_{+}}}

\newcommand{\borel}[1]{\mathcal{B} (#1 )}

\newcommand{\BigO}[1]{\mathop{}\!O{\left(#1\right)}}

\newcommand{\predictableVariation}[1]{ \langle#1\rangle }

\newcommand{\absolute}[1]{| #1 | }

\makeatletter
\let\oldabs\abs
\def\abs{\@ifstar{\oldabs}{\oldabs*}}

\newcommand{\cadlag}{c\`adl\`ag }

\newcommand{\metricGivenSpace}[1]{\mathsf{d}_{#1} }

\newcommand{\defeq}{\coloneqq}

\newcommand{\indicator}[2]{\mathsf{1}_{#1}{\left(#2\right)}}

\newcommand{\differential}[1]{\mathrm{d} #1}
\newcommand{\timeDerivative}[1]{\frac{\differential}{\differential t} #1 }

\newcommand{\effectiveDomain}[1]{\mathsf{dom}\, #1}

\newcommand{\eqstop}{.}
\newcommand{\eqcomma}{,}

\newcommand{\norm}[1]{\left\lVert#1\right\rVert}

\newcommand{\E}{\mathsf{E}}
\newcommand{\Eof}[1]{\E\left[#1 \right]}

\newcommand{\prob}{\mathsf{P}}
\newcommand{\probOf}[1]{\prob\left(#1\right)}

\newcommand{\history}[1]{\mathcal{F}_{#1}  }

\newcommand{\disteq}{\, \overset{\mathcal{D}}= \, }

\newcommand{\ConvInProb}{\xrightarrow[]{ \hspace*{4pt}  \text{         P   } }}

\newcommand{\ConvInDist}{  \overset{\hspace*{4pt}  \mathcal{D}  }{\implies } }

\newcommand{\myExp}[1]{\exp \left( #1 \right)  }

\newcommand{\ie}{\textit{i.e.}}

\newcommand{\nX}{X^{(n)}}
\newcommand{\nY}{Y^{(n)}}

\title{Stochastic Analysis of Entanglement-assisted Quantum Communication Channels} 

\author{Karim S. Elsayed}
\address{Karim~S.~Elsayed, 
Institute of Communications Technology, Leibniz Universit\"at Hannover, Appelstra{\ss}e 9a,
30167 Hanover, Germany
}
\email{karim.elsayed@ikt.uni-hannover.de}

\author{Olga Izyumtseva}
\address{Olga Izyumtseva, 
School of Mathematical Sciences, 
University of Nottingham, 
University Park, 
Nottingham NG7 2RD, 
United Kingdom
}
\email{Olga.Iziumtseva1@nottingham.ac.uk}

\author{Wasiur R. KhudaBukhsh}
\address{Wasiur R. KhudaBukhsh, 
School of Mathematical Sciences, 
University of Nottingham, 
University Park, 
Nottingham NG7 2RD, 
United Kingdom
}
\email{wasiur.khudabukhsh@nottingham.ac.uk}

\author{Amr Rizk}
\address{Amr Rizk, 
Institute of Communications Technology, Leibniz Universit\"at Hannover, Appelstra{\ss}e 9a,
30167 Hanover, Germany
}
\email{amr.rizk@ikt.uni-hannover.de}

\begin{document}
\maketitle

\begin{abstract}
    We present a queueing model for a quantum communication network consisting 
    of a primary queue and a service queue in which Bell pairs are formed and stored. The Bell pairs are inherently extremely short-lived rendering the service queue (the quantum queue) much faster than the primary queue. We study the asymptotic behaviour of this multi-scale queueing system via a  stochastic averaging principle. We prove a \ac{FLLN} and a \ac{FCLT} for the standard queue averaging the dynamics of the fast service queue. 
\end{abstract}



\tableofcontents

\section{Introduction}
    \label{sec:intro}
Recently, quantum communication has gained significant interest stemming from its potential to improve current distributed computing applications beyond their classical limits or enable quantum-specific applications.
Quantum communication describes the transmission of quantum particles, called qubits, through quantum channels, e.g., the transmission of photons through an optical fibre~\cite{cacciapuoti_Caleffi_entang_meets_calssical,bennett1999entanglement}. The unique resources in quantum communication responsible for its potential capabilities are called entanglements. 
An entanglement in its simplest form is a correlation shared between two qubits. 
It is known that utilizing an entanglement shared between two communicating nodes may increase the classical capacity of the quantum channel between them to a maximum of double the classical channel capacity~\cite{bennett1999entanglement,bennett2002entanglement,hao2021entanglement}.
In this paper, we consider a communication system that takes advantage of sporadic entanglements to assist classical communication. 
Here, messages of random size are transmitted using entanglement assistance whenever possible, otherwise, they are transmitted classically. 
A rigorous asymptotic analysis of such communication systems, a rapidly growing field driven by recent experimental successes~\cite{yuan23quantumrouter,azuma2023quantum,fustionbasedqcomp}, is missing in the literature.

\begin{figure}[h]
    \centering
    \includegraphics[ trim={0cm 1cm 0cm 0cm},clip]{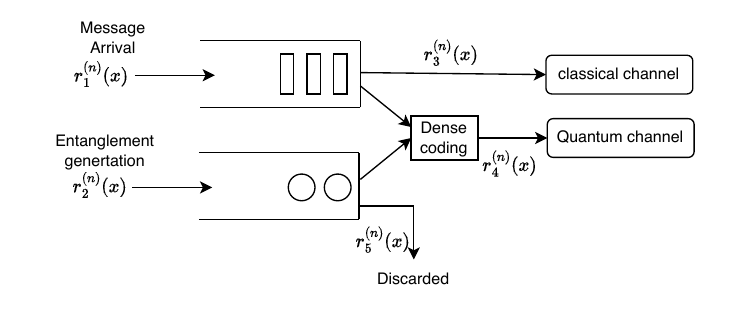}
    \caption{Queueing system: messages are sent through the quantum channel with an accelerated service propensity function $r_4^{(n)}(x)$ only when entanglements (in queue B) are available, otherwise messages are served with rate $r_3^{(n)}(x)$ through the classical communication channel. Queue B gets depleted at rate $r_5^{(n)}(x)$ due to the quality decay of entanglements with time. Here, $x=(x_1,x_2)$ is the vector of queue lengths of the messages and entanglement queues and the functions $r_{1}^{(n)}(x), r_{2}^{(n)}(x), r_{3}^{(n)}(x), r_{4}^{(n)}(x),$ and  $r_{5}^{(n)}(x)$ are the propensity/intensity functions defined in \Cref{eq:propensity_1} and \Cref{eq:propensity_2} in \Cref{sec:ctmc_model}. The superscripts on the rate functions are included to emphasise their dependence on some scaling parameter $n$. 
    }
    \label{fig:queueing_system}
\end{figure}

We present a queueing model 
consisting of two queues as depicted in~\Cref{fig:queueing_system}. The queue lengths are described by a system of  \acp{SDE} driven by \acp{PRM}. 
Queue A (the primary queue) stores the information messages with an arrival rate  $r_1(.)$ and queue B stores the entanglements, which are created and shared between the two nodes according to the rate function $r_2(.)$.
We define all rate (sometimes called propensity) functions precisely in \Cref{sec:ctmc_model}, and provide a discussion about their practical implementation in \Cref{sec:discussion}.
When available, entanglements encode the classical information through a class of dense coding~\cite{bose2000mixed}, offering a service speed up according to a service propensity function $r_4(.)$.
Otherwise, messages are transmitted through the classical channel according to classical service rate $r_3(.)$. 
In addition, entanglements are depleted from Queue~B at rate $r_5(.)$ since the usability of entanglements is known to be limited to their small lifetimes due to the impact of quantum noise~\cite{schlosshauer_decoherence}.
Note that our analysis stands in contrast to classical capacity calculation of quantum channels~\cite{nielsen_QC_QI_book} since we do not consider the limiting behaviour of the system in terms of the utilisation of the channel capacity. 

\paragraph{Mathematical contributions}
The literature on the stochastic analysis of quantum communication channels is sparse. This is due to the recency of the technology but also due to the nontriviality of the required translation between the quantum terminology and the language of classical probability theory. 
Our primary contribution in this paper is an asymptotic analysis of the scaled queue length at the primary queue based on the theory of stochastic averaging principle \cite{Kurtz1992Averaging,Hasminskii1966a,Khasminskii1966b}. To be precise, we prove a \acf{FLLN} (\Cref{thm:relative_compactness}) and a \acf{FCLT} (\Cref{thm:clt}) for the scaled queue length at the queue A when the dynamics of the faster service queue (the queue B) are averaged at the timescale of the standard queue. 

Our proofs are direct and probabilistic in nature as opposed to the often used analytic generator-based methods. 
We rely on the stochastic analysis of the paths of \acp{SDE} driven by \acp{PRM} \cite{Applebaum_2009Levy}, and various martingales associated with them. The main technical challenge in proving the \ac{FLLN} lies in verifying the tightness of the occupation measure of the fast process, and identifying the unique limiting measure. For the \ac{FCLT}, we explicitly solve a Poisson equation associated with the generator of the fast process (with the slow variable frozen) and use the \ac{MCLT} along with the tightness-uniqueness arguments for weak convergence. The covariation process of the limiting Gaussian semimartingale is obtained explicitly thanks to the explicit solution of the Poisson equation.

\paragraph{Connection to existing literature on stochastic averaging principle}
The ideas of stochastic averaging date back to Khas'minskii \cite{Hasminskii1966a,Khasminskii1966b} with many important follow-up applications to stochastic approximation algorithms~\cite{Kushner2009stochApprox}, queueing networks \citep{Perry2013Fluid}, and stochastic \acp{CRN} \cite{Kang:2013:STM,Kang:2014:CLT,Ball:2006:AAM,Enger2023Unified,Crudu2012AAP,ganguly2024enzymekinetic,Ganguly2026tQSSA,TK2025averaging,fromion2025stochastic,hashemi2016stochastic}. The standard road map to proving a stochastic averaging principle has been known in the probability theory community for several decades. Notably, the late probabilist Thomas G. Kurtz provided a general recipe for proving stochastic averaging principle for martingale problems in \cite{Kurtz1992Averaging}. The proof of our \ac{FLLN} (\cref{thm:relative_compactness}) is based on this approach of Kurtz. The main result of \cite{Kurtz1992Averaging} is restated as \cref{thm:stochastic_averaging} in \cref{sec:math_details} for the readers' convenience. The first rigorous application of the stochastic averaging principle to \acp{CRN} appears in \cite{Ball:2006:AAM}, where many specific examples including a stochastic model of intracellular viral kinetics are studied. Later Kang and Kurtz provide a systematic approach to proving \acp{FLLN} under time-scale separation by verifying a number of ``balance conditions''  for general \acp{CRN} in \cite{Kang:2013:STM}. The corresponding \acp{FCLT}, akin to our \Cref{thm:clt}, were discussed in \cite{Kang:2014:CLT}. There, the authors list a number of sufficient conditions for the convergence of an appropriately scaled slow component to a Gaussian semimartingale. The main result of \cite{Kang:2014:CLT} (Theorem 2.11 in the paper) is stated in terms of the solution of a Poisson equation, as is ours too. However, we solve the required Poisson equation explicitly, which is often not feasible and makes our result immediately applicable.  Note that verifying the sufficient conditions in \cite{Kang:2013:STM,Kang:2014:CLT,Ball:2006:AAM} for our queueing model is not straightforward. For example, a rigorous validation of \cite[Condition 2.1 - Condition 2.10]{Kang:2014:CLT} requires a detailed, and almost a paper-length analysis.  Similar to the current paper, much of \cite{Kang:2013:STM,Kang:2014:CLT} is based on Kurtz's previous work \cite{Kurtz1992Averaging}, and an efficient use of the \ac{MCLT} along with the solution to a Poisson equation.

Besides the works of Kurtz and colleagues, there have been several works that propose multiscale approximations utilising the stochastic averaging principle (sometimes with some overlap with Kurtz's approach). For instance, the authors in \cite{Enger2023Unified} propose a generator-based weak convergence theorem, which incidentally is also based on \cite[Theorem 2.1]{Kurtz1992Averaging}. Our result does not follow from \cite[Theorem 2.3]{Enger2023Unified} since the operator $\mathcal{A}$ (defined in \eqref{eq:A_operator}) describing the dynamics of the slow variable depends on both the fast and the slow variables. The work of \cite{Crudu2012AAP} provides stochastic averaging principles where the limiting process is a piecewise deterministic Markov jump process. The works \cite{Laurence2025AIMD,Laurence2026ScalingMethods}  consider \acp{CRN} with non-standard scaling parameters such as the initial volume of the \acp{CRN} and discontinuous limits.


To the best of our knowledge, this is the first application of the stochastic averaging principle to entanglement-assisted communication networks. 
The results of this paper provide the necessary ingredients for the analysis of more complex entanglement queueing models for quantum communication networks.
One promising direction is the analysis of multi-user (or multi-receiver) scenarios, where either users share the same entanglement queue or each user has its own dedicated entanglement queue, as recently developed in~\cite{Dai_Qu_Queuing_delay}. 
Similar multi-entanglement queueing models arise in the context of quantum switches, the counterpart of classical switches, and are currently being developed in several recent works~\cite{Towsley_stochastic_qu_switch,Twosley_ideal_Qu_switch}.
Among the recent works, \cite{bhambay2025optimal,zubeldia2026matching} rely on similar principles, namely time scale separation and~\ac{FLLN}, however, in the context of scheduling in multi-entanglement queueing scenarios.
%
To help facilitate practical applications, we have included an extensive discussion of various model parameters and their implications in \Cref{sec:discussion}. For this reason, we have also made an effort to keep the mathematical exposition as self-contained as possible.


    The rest of the paper is structured as follows. The \ac{SDE} model is described in \Cref{sec:ctmc_model}. The statement and the proof of the \ac{FLLN} is provided in \Cref{sec:averaging}. We discuss the \ac{FCLT} in \Cref{sec:FCLT}.  Finally, we conclude with a discussion in \Cref{sec:discussion}. Additional mathematical details are provided in \ref{sec:math_details}.  Before delving into the mathematical model, we lay down our notational conventions. 

    \paragraph{Notational conventions} The indicator (characteristic) function of a set $E$ is denoted by $\indicator{E}{\cdot }$, \ie, $ \indicator{E}{x} =1$  if  $x\in E$, and zero otherwise. 
    Given a complete, separable metric space $(E, \metricGivenSpace{E})$, we will denote by $\mathcal{M}(E)$ the space of finite measures on $E$ with the weak topology. Given an interval $E \subseteq [0, \infty)$ and a metric space $F$, 
    we will use $C( E, F)$, and $D(E, F) $ to denote the space of continuous $F$-valued functions on $E$, and the space of $F$-valued functions on $E$ with \cadlag paths, respectively. We will endow the space $C$ with the topology of the uniform norm and the space $D$ with the Skorohod topology \cite{Billingsley1999Convergence,jacod2003limit}. When $F$ is a subset of an Euclidean space, the subspace of $C(E, F)$ containing bounded and continuous functions will be denoted by $C_b(E, F)$.  Furthermore, the space of infinitely differentiable (smooth) $F$-valued functions of $E$ will be denoted as $C^{(\infty)}(E, F)$, with the subspace $C_{c}^{(\infty)}(E, F) \subset C^{(\infty)}(E, F)$ containing those with compact supports. 
    Let us denote by  ${\mathcal{M}}_1([0, \infty)\times E)$ the space of non-negative measures on $[0, \infty)\times E$ with the property that $\mu([0, t]\times E ) =t $ for every $t\ge 0$ and every $\mu \in {\mathcal{M}}_1([0, \infty)\times E)$. Denote the Prohorov metric on $\mathcal{M}([0, t]\times E)$ by $q_t$ \cite{Billingsley1999Convergence,Ethier:1986:MPC}. Then, the space $\mathcal{M}_1([0, \infty)\times E)$ will be equipped with the metric 
    \begin{align*}
        q(\mu, \nu) \defeq \int_{0}^{\infty} e^{-t} \min\{1, q_t(\mu, \nu) \}\differential{t} \eqcomma
    \end{align*} 
    where we apply $q_t$ on the restrictions of $\mu$ and $\nu$ to $[0, t]\times E$. 
    

\section{Stochastic Model}
    \label{sec:ctmc_model}
    Consider the entanglement-assisted queueing system depicted in \Cref{fig:queueing_system}. Throughout this paper, we will use the letter $n$ as a scaling parameter, and the limiting results will be understood as $n\to \infty$.   
     The scaling factor $n$ reflects the separation in time scales between the entanglement generation at queue $B$ and message arrivals at queue $A$, which differ by several orders of magnitude. 
    This arises from the fast generation of entanglements between two nodes, typically involving the transmission of a single photon, in contrast to classical messages that evolve on longer time scales governed by message sizes, network congestion and latency. The use of such a non-physical asymptotic large parameter $n$ (or a small parameter $1/n$) is common in the literature on fluid and diffusion limits of queues \citep{Perry2013Fluid}, multiscale approximation in \acp{CRN} \cite{Kang:2013:STM,Ganguly2026tQSSA,Kang:2014:CLT}.  
     
    Let $\nX (t) \defeq (X_A^{(n)}(t), X_B^{(n)}(t))$, where $X_A^{(n)}(t)$, and $X_B^{(n)}(t)$ are the queue lengths at the queues $A$ and~$B$ respectively at time $t$. Assume that $\nX (0)=(0, 0)$, \ie, the queues are empty initially. This assumption is made only for the sake of simplicity, and does not affect the generality of our results as long as the initial state is independent of parameter $n$. If it does depend on $n$, we need to make additional assumptions about how the initial state scales with $n$. Since it does not bring any new mathematical challenges, we simply assume $\nX(0)=(0,0)$. 
    In order to describe the dynamics of the process $\nX$, let us define the rate functions 
    \begin{align}
    \label{eq:propensity_1}
        r_1^{(n)}(x) \defeq nr_1\left(\frac{x_1}{n}\right),\;\quad r_2^{(n)}(x) \defeq n\lambda,\; \quad r_3^{(n)}(x) \defeq  nr_3\left(\frac{x_1}{n}\right) \indicator{\{0\}}{x_2},\; 
    \end{align}
    and 
    \begin{align}
    \label{eq:propensity_2}
        r_4^{(n)}(x) \defeq n r_{4}\left(\frac{x_1}{n}\right)x_2 
        , \; \quad r_5^{(n)}(x) \defeq n\mu x_2
        \eqcomma 
    \end{align}
    for $x\defeq (x_1, x_2) \in \mathbb{N}^2_{0}$, $\lambda>0,\mu>0,$ and some functions $r_1,\  r_3,\  r_{4}:\mathbb{R}_{+}\to\mathbb{R}_{+}$, which are assumed Lipschitz continuous. In order to avoid trivialities and absorbing states, we assume $r_1(0)>0$, and $r_3(0)= r_4(0)=0$. The choice of $r_2^{(n)}$ and $r_5^{(n)}$ reflects the fact that the fast queue would behave like an $M/M/\infty$ queue in the absence of queue $A$.
    Recall that the rate functions $r_1(.)$ and $r_2(.)$ model the message arrival in queue A and the entanglement generation in Queue B. 
    Moreover, the rate functions $r_3(.)$ and $r_4(.)$ describe the classical service and the entanglement assisted service, respectively, while the function $r_5(.)$ models the (fast) depletion of entanglements from queue B due to quantum noise.  
%

    
    The process $X^{(n)} $ has paths in  $D([0,\infty), \mathbb{N}^2_{0})$, the space of \cadlag functions on $[0, \infty)$ endowed with the Skorohod topology \cite{Ethier:1986:MPC,Billingsley1999Convergence}.
    We describe the trajectories of the stochastic process $X^{(n)}$ by means of the following \acp{SDE} \cite{IkedaWatanabe2014Stochastic,Applebaum_2009Levy}     
    \begin{align}
        \begin{aligned}
            X_A^{(n)}(t) &{} = X_A^{(n)}(0) + \int_{0}^{t} \int_{0}^{\infty}\indicator{[0, \; r_1^{(n)} (\nX(u-)) ]}{v} Q_{A}^{(1)}(\differential{u}, \differential{v}) \\
            &{} \quad - \int_{0}^{t} \int_{0}^{\infty}\indicator{[0, \; r_3^{(n)} (\nX(u-)) ] }{v} Q_A^{(2)}(\differential{u}, \differential{v})  \\
            &{} \quad - \int_0^t \int_0^\infty \indicator{[0, \; r_4^{(n)} (\nX(u-)) ]}{v}Q_{AB}^{(1)}(\differential{u}, \differential{v}) \eqcomma  \\ 
            X_B^{(n)}(t) &{} = X_B^{(n)}(0) + \int_{0}^{t} \int_{0}^{\infty}\indicator{[0, \; r_2^{(n)} (\nX(u-))]}{v} Q_{B}^{(1)}(\differential{u}, \differential{v}) \\
            &{} \quad 
            - \int_{0}^{t} \int_{0}^{\infty}\indicator{[0, \;  r_5^{(n)} (\nX(u-)) ]}{v} Q_B^{(2)}(\differential{u}, \differential{v}) \\
            &{} \quad 
            - \int_0^t \int_0^\infty \indicator{[0, \; r_4^{(n)} (\nX(u-)) ]}{v}Q_{AB}^{(1)}(\differential{u}, \differential{v}) \eqcomma 
        \end{aligned}
        \label{eq:nX_sde}
    \end{align}
   where $Q_A^{(1)}(\differential{u}, \differential{v}), Q_{A}^{(2)}(\differential{u}, \differential{v})$, $Q_B^{(1)}(\differential{u}, \differential{v}), Q_{B}^{(2)}(\differential{u}, \differential{v})$, and $Q_{AB}^{(1)}(\differential{u}, \differential{v})$  are independent \acp{PRM} on $\setOfPositiveReals\times \setOfPositiveReals$ with intensity $\differential{u}\times \differential{v}$ where $\differential{u}$, and $\differential{v}$ are Lebesgue measures on $\setOfPositiveReals$. The random measures $Q_A^{(1)}, Q_A^{(2)}, Q_B^{(1)}$, $Q_B^{(2)}$, and $Q_{AB}^{(1)}$ are defined on the same probability space $(\Omega, \history{}, \prob)$, and are independent of $X^{(n)}(0)$. Assume $\history{}$ is $\prob$-complete and  associate to $(\Omega, \history{}, \prob)$ the filtration $(\history{t}^n)_{t\ge 0}$ generated by $\nX(0), Q_{A}^{(1)}, Q_A^{(2)}, Q_B^{(1)}, Q_B^{(2)}$, and $Q_{AB}^{(1)}$.
    Let $\history{0}$ contain all $\prob$-null sets in $\history{}$. 

    \section{\acl{FLLN}}
    \label{sec:averaging}
    We will investigate the scaling behaviour of the slow queue $A$ by averaging the fast queue $B$. Therefore, let us define the scaled stochastic process 
    \begin{align*}
        {\nY} = (\nY_A, \nY_B) \defeq (n^{-1} \nX_A, \nX_B)\eqstop 
    \end{align*}
    Then, the stochastic process $\nY$ is a continuous time Markov process with the generator (\cite{Ethier:1986:MPC}, see also \cite[Appendix 1]{Kipnis1999Scaling})
    \begin{align}
        \begin{aligned}
        \mathcal{L}_n f(y) & \defeq  n r_1(y_1)\left(f(y_1 +\frac{1}{n}, y_2) - f(y)\right) +  n \lambda\left(f(y_1, y_2+1) - f(y)\right) \\
        &{}\quad  
        + n r_3(y_1) \indicator{\{0\}}{y_2} \left(f(y_1-\frac{1}{n}, y_2) - f(y)\right) 
        \\
        &{}\quad 
        + n r_{4}(y_1) y_2  \left(f(y_1-\frac{1}{n}, y_2-1) - f(y)\right) 
        + n \mu y_2 (f(y_1, y_2-1) - f(y)), 
        \end{aligned}
           \label{eq:generator}
    \end{align}
    for bounded continuous functions $f: \ \mathbb{R}_{+}\times\mathbb{N}_0 \to \setOfReals$, and $y \defeq (y_1, y_2) \in\mathbb{R}_{+}\times\mathbb{N}_0$. Moreover, for $f\in C_b(\mathbb{R}_{+}\times\mathbb{N}_0,\setOfReals)$, the stochastic process $M_f^{(n)}$ defined by 
    \begin{align}
        M_f^{(n)}(t) \defeq f(\nY(t)) - f(\nY(0)) - \int_{0}^{t} \mathcal{L}_n f(\nY(s))\differential{s}\, \eqcomma \quad t\ge 0 \eqcomma 
        \label{eq:scaled_martingale}
    \end{align} 
    is an $\history{t}^n$-martingale \cite{Ethier:1986:MPC,Kipnis1999Scaling}, called the Dynkin's martingale. In other words, the Markov process $\nY$ solves the martingale problem for $(\mathcal{L}_n, \effectiveDomain{\mathcal{L}_n})$ with $\effectiveDomain{\mathcal{L}_n} \defeq C_b(\mathbb{R}_{+}\times\mathbb{N}_0,\setOfReals)$. Moreover, it follows from \Cref{eq:nX_sde} that the stochastic process $\nY$ evolves according to the following \acp{SDE}
    \begin{align}
        \begin{aligned}
            \nY(t) &{} = \nY(0) + \frac{1}{n}\int_{0}^{t} \int_{0}^{\infty}\indicator{[0, \; n r_1(\nY_A(u-)) ]}{v} Q_{A}^{(1)}(\differential{u}, \differential{v}) \\
            &{} \quad - \frac{1}{n} \int_{0}^{t} \int_{0}^{\infty}\indicator{[0, \; n r_3 (\nY_A(u-)) \indicator{\{0\}}{\nY_B(u-)} ] }{v} Q_A^{(2)}(\differential{u}, \differential{v})  \\
            &{} \quad - \frac{1}{n} \int_0^t \int_0^\infty \indicator{[0, \; n r_4 (\nY_A(u-))\nY_B(u-) ]}{v}Q_{AB}^{(1)}(\differential{u}, \differential{v}) \eqcomma  \\ 
            \nY_B(t) &{} = \nY_B(0) + \int_{0}^{t} \int_{0}^{\infty}\indicator{[0, \; n \lambda]}{v} Q_{B}^{(1)}(\differential{u}, \differential{v}) \\
            &{} \quad 
            - \int_{0}^{t} \int_{0}^{\infty}\indicator{[0, \;  n \mu \nY_B(u-) ]}{v} Q_B^{(2)}(\differential{u}, \differential{v}) \\
            &{} \quad 
            - \int_0^t \int_0^\infty \indicator{[0, \; n r_4 (\nY_A(u-))\nY_B(u-) ]}{v}Q_{AB}^{(1)}(\differential{u}, \differential{v}) \eqstop 
        \end{aligned}
        \label{eq:nY_sde}
    \end{align}
    
    The form of the generator $\mathcal{L}_n$ in \Cref{eq:generator} and the \acp{SDE} in  \Cref{eq:nY_sde} reveal that the stochastic process $\nY_B$ jumps rapidly (at rate $\BigO{n}$), while the stochastic process $\nY_A$ has approximately deterministic dynamics driven by the \emph{average} dynamics of $\nY_B$. To describe the dynamics of the fast process, let us define the linear operator 
    \begin{align}
        \mathcal{B}_{y_1}h(y_2) \defeq \lambda(h(y_2+1)-h(y_2)) + (r_{4}(y_1)  + \mu )y_2 (h(y_2-1) - h(y_2)) \eqcomma 
       \label{eq:frozen_generator}
    \end{align}
    for each fixed $y_1$, and  for bounded functions $h: \setOfNonnegativeIntegers \mapsto \setOfReals$. It is straightforward to verify that the linear  operator $\mathcal{B}_{y_1}$ generates a pure jump Markov process, namely, a \ac{BD} process with birth rate $\lambda$, and death rate $(r_{4}(y_1) + \mu )y_2,\ y_2\geqslant 0.$ 
       Furthermore, for all $y_1 \ge 0$,  it is ergodic
    and admits a unique stationary distribution $\pi(\cdot)$ \cite{Robert2003StochasticNetworks} given by 
    \begin{align}
        \begin{aligned}
            \pi(k) & = \myExp{-m(y_1)} \frac{m(y_1)^k }{k!}\eqcomma  \text{ for } k\ge 0\eqcomma
            \text{with } m(y_1)  = \frac{\lambda}{(r_{4}(y_1)  + \mu)}\eqstop 
        \end{aligned}
        \label{eq:stationary_distribution}
    \end{align}
    Interestingly, despite its rapid jumps, it takes a long time for the process $\nY_B$ to reach very large queue lengths. 
    More precisely, define the hitting time random variables   
    \begin{align}
        \tau^{(n)}_k&{} \defeq \inf\{t\geqslant 0:\ Y^{(n)}_B(t) = k\}\eqcomma 
        \label{eq:hitting_times}
    \end{align}
    for $k \ge 1$. The following lemma gives an estimate on the \ac{MGF} of the hitting times from below.
    \begin{myLemma}
        Assume $\nY_B(0) = 0$, and set $\rho \defeq \lambda/\mu$.  Then, for $\beta > 0$, and $k \ge 1$,  we have the following estimate 
        \begin{align*}
            \Eof{\myExp{-\beta \tau^{(n)}_k}} \le \bar{\phi}_{\beta}(n, k) \defeq  \frac{\int^{\infty}_0u^{\frac{\beta}{n\mu}-1}\myExp{-\rho u}\differential{u}}{\int^{\infty}_0(1+u)^k u^{\frac{\beta}{n\mu}-1}\myExp{-\rho u}\differential{u}} \eqstop 
        \end{align*}
        \label{lem:MGF_hitting_times}
    \end{myLemma}
    \begin{proof}[Proof of \Cref{lem:MGF_hitting_times}]
        Let us define the functions
 $$
 g^{(n)}_{\alpha}(t,y)=(1+\alpha \myExp{n\mu t})^{y_2}\myExp{-\rho\alpha\myExp{n\mu t}},
 $$
for $\alpha\in\mathbb{R},\ n\geqslant 1,$ and $y=(y_1,y_2)\in \setOfPositiveReals\times\setOfNonnegativeIntegers.$ Then, note that 
\begin{align*}
\mathcal{L}_ng^{(n)}_{\alpha}(t,\cdot)(y)
&{}=e^{-\rho\alpha\exp{(n\mu t)}}n\alpha e^{n\mu t}\left[\lambda(1+\alpha e^{n\mu t})^{y_2} 
-(\mu+r_4(y_1))y_2(1+\alpha e^{n\mu t})^{y_2-1}\right].
\end{align*}
Furthermore, we have 
\begin{align*}
\frac{\partial}{\partial t}g^{(n)}_{\alpha}(t,y)&{}=e^{-\rho\alpha\exp{(n\mu t)}}n\alpha e^{n\mu t}\left[\mu y_2(1+\alpha e^{n\mu t})^{y_2-1}-\lambda(1+\alpha e^{n\mu t})^{y_2}\right]\\
&{}\leqslant e^{-\rho\alpha\exp{(n\mu t)}}n\alpha\left[(\mu+r_4(y_1)) y_2(1+\alpha e^{n\mu t})^{y_2-1}-\lambda(1+\alpha e^{n\mu t})^{y_2}\right]\\
&{}\leqslant-\mathcal{L}_ng^{(n)}_{\alpha}(t,\cdot)(y).
\end{align*}
Hence, for any $(t,y)\in\setOfPositiveReals \times (\setOfPositiveReals\times \setOfNonnegativeIntegers)$, we have 
$$
\frac{\partial}{\partial t}g^{(n)}_{\alpha}(t,y)+\mathcal{L}_ng^{(n)}_{\alpha}(t,\cdot)(y)\leqslant 0,
$$
\ie,  the functions $g^{(n)}_\alpha$ are space-time superharmonic for the generator $\mathcal{L}_n$ (see \Cref{def:space_time_subsuperharmonic} in Appendix~\ref{sec:math_details}). 
Then, by virtue of \Cref{lemma:harmonic_martingale} in Appendix~\ref{sec:math_details},  the  process
$$
M^{(n)}_g(t)\defeq g^{(n)}_{\alpha}(t,Y^{(n)}(t))\eqcomma\quad  t\ge 0\eqcomma 
$$
is a local supermartingale. In fact, it is a supermartingale for each fixed $n$, since
$$
Y^{(n)}_B(t)\leqslant\int^t_0\int^{\infty}_0\indicator{[0,n\lambda]}{v}Q^{(1)}_B(du,dv)\disteq N(n\lambda t)
$$
for a unit rate Poisson process $N.$ Therefore, for any $t\geqslant 0$, we have 
$$
\Eof{M^{(n)}_g(t)}\leqslant g^{(n)}_{\alpha}(0,Y^{(n)}(0))=e^{-\rho\alpha}\eqcomma 
$$ 
which implies that
$
\Eof{(1+\alpha \myExp{n\mu t})^{Y^{(n)}_B(t)}}\leqslant \myExp{\rho\alpha(\myExp{n\mu t}-1)} \eqstop 
$
By a change of variable $\theta=1+\alpha \myExp{n\mu t}$, one gets
$
\Eof{\theta^{Y^{(n)}_B}}\leqslant \myExp{\rho(\theta-1)(1-\myExp{-n\mu t})}.
$
Using the above estimate and following arguments similar to those used to prove \cite[Proposition 6.6, p. 148]{Robert2003StochasticNetworks}, we can verify that for each $\alpha>0,\ t\in\mathbb{R}_{+},$ the  process
$$
\phi^{(n)}_\alpha(t)=e^{-\alpha n\mu t}\int^{\infty}_0(1+u)^{Y^{(n)}_B(t)}u^{\alpha-1}e^{-\rho u}\differential{u}\eqcomma \quad t \ge 0\eqcomma 
$$
is a supermartingale. Therefore, for the stopping times $\tau^{(n)}_k$ defined in \Cref{eq:hitting_times},  we immediately have the following estimate 
$$
\Eof{\phi^{(n)}_\alpha(\tau^{(n)}_k)}\leqslant\Eof{\phi^{(n)}_\alpha(0)}=\int^{\infty}_0u^{\alpha-1}e^{-\rho u}\differential{u}.
$$
However, we know $Y^{(n)}_B(\tau^{(n)}_k)=k$. Therefore,  we get 
\begin{align*}
\Eof{\myExp{-\alpha n\mu\tau^{(n)}_k}\int^{\infty}_0(1+u)^k u^{\alpha-1}e^{-\rho u}\differential{u} }\leqslant\int^{\infty}_0u^{\alpha-1}e^{-\rho u}\differential{u}.
\end{align*}
Therefore, for  $\beta>0$, we have 
\begin{align*}
\Eof{e^{-\beta\tau^{(n)}_k}}&{}=\Eof{e^{-\left(\frac{\beta}{n\mu}\right)n\mu\tau^{(n)}_k}}
\leqslant\frac{\int^{\infty}_0u^{\frac{\beta}{n\mu}-1}\myExp{-\rho u}\differential{u}}{\int^{\infty}_0(1+u)^k u^{\frac{\beta}{n\mu}-1}\myExp{-\rho u}\differential{u}} = \bar{\phi}_{\beta}(n,k)\eqcomma 
\end{align*}
which completes the proof. 
    \end{proof}

    \begin{myRemark}
        Notice that for each fixed $n$, due to the \ac{MCT}, we have 
$$
\lim_{k\to\infty}\int^{\infty}_0(1+u)^k u^{\frac{\beta}{n\mu}-1}\myExp{-\rho u}\differential{u}=\infty.
$$
Hence, for each fixed $n$,  
$
\lim_{k\to\infty}\bar{\phi}(n,k)=0.
$
Moreover, for each $k\geqslant 1,$ again due to the \ac{MCT}, we have
$$
\lim_{n\to\infty}\bar{\phi}(n,k)=\frac{\int^{\infty}_0u^{-1}\myExp{-\rho u}\differential{u}}{\int^{\infty}_0(1+u)^k u^{-1}\myExp{-\rho u}\differential{u}} \in [0,1]\eqstop 
$$
\label{remark:MGF_hitting_times}
    \end{myRemark}
    
    \Cref{lem:MGF_hitting_times} implies that it takes a long time for the process $\nY_B$ to reach very large queue lengths, which is useful to verify a compact containment condition. However, because of the rapid jumps of $\nY_B$, the process $\nY_B$ fails to satisfy the modulus of continuity condition (Condition (ii) in \cite[Chapter 3, Theorem 13.2]{Billingsley1999Convergence}) and hence, is not relatively compact as a $D([0, \infty), \setOfNonnegativeIntegers)$-valued random element. However, its occupation measure is a better behaved quantity and is indeed relatively compact as a collection of random measures. We will see that the occupation measure converges to a deterministic measure, which depends on the stationary distribution corresponding to the generator $\mathcal{B}_{y_1}$ as $n\to \infty$ when the slow variable is frozen at $y_1$. Therefore, the dynamics of $\nY_B$ can be approximated by its ergodic moments. We will now make the above heuristic argument more precise in the following.


    Consider the occupation measure $\Gamma_n$ of the fast variable $\nY_B$, defined as 
    \begin{align*}
        \Gamma_n ([0, t]\times F) \defeq \int_{0}^{t} \indicator{F}{\nY_B(s) }\differential{s} \eqcomma 
    \end{align*} 
    for $t>0$, and $F \in \mathcal{B}(\setOfNonnegativeIntegers) \defeq 2^\setOfNonnegativeIntegers$. Naturally, $\Gamma_n$ is a random measure, which we view as an $\mathcal{M}_1([0, \infty)\times\mathbb{N}_0)$-valued random element. We have the following \ac{FLLN}. 
    
    

    \begin{myTheorem}[\acl{FLLN}]
        Assume $\nY(0) = (0,0)$, and $T\in (0, \infty)$. Then, the sequence $\{\nY_A : n\ge 1\}$ converges in probability to a deterministic function $y_A$ in $D([0, T], \setOfPositiveReals)$, where the function $y_A$ lies in $C([0, T], \setOfPositiveReals)$ almost surely, and  evolves according to  the following \ac{ODE}
   \begin{align}
   \timeDerivative{y_A}
   = G(y_A) \defeq  r_1(y_A) -r_3(y_A) \myExp{-m(y_A)} -r_4(y_A)m(y_A)\eqcomma 
   \label{eq:limit_ODE}
   \end{align}
   with the initial condition $y_A(0)=0$.
       \label{thm:relative_compactness}
    \end{myTheorem}


     \begin{proof}[Proof of \Cref{thm:relative_compactness}]  
        We will proceed in two steps. We first show that
 the  sequence $\{(\nY_A, \Gamma_n) : n\ge 1\}$ is relatively compact in $D([0, T], \setOfPositiveReals)\times \mathcal{M}_1([0, T]\times \mathbb{N}_{0}) $, and then identify the limit. 

 \paragraph{Relative compactness}
Note that from \Cref{eq:nY_sde}, 
         \begin{align*}
            \sup_{t\in[0,T]}\nY_A(t) &{} \leqslant \tilde{Y}^{(n)}_A(T) = \frac{1}{n}\int_{0}^{T} \int_{0}^{\infty}\indicator{[0, \; r^{(n)}_1(nY^{(n)}_A(u-),Y^{(n)}_B(u-))]}{v} Q_{A}^{(1)}(\differential{u}, \differential{v}),
                         \end{align*}
       since  $ \nY_A(0)=0$, and $\tilde{Y}^{(n)}_A$ is a non-decreasing process.         
        To check the compact containment condition (see \Cref{def:compact_containment}, and \cite{Ethier:1986:MPC}) for $\nY_A$,
        %
        let us show that for each $\epsilon>0$ there exists a positive constant $k_{\epsilon}$ such that
        \begin{align}
            \inf_{n}\probOf
        {\sup_{t\in[0,T]}\nY_A(t)\leqslant k_{\epsilon}}>1-\epsilon. 
        \nonumber
        \end{align}
        Note that 
        \begin{align*}
            \Eof{\nY_A(t)} \leqslant \int^T_0 \Eof{r_1(Y^{(n)}_A(t))}\differential{t}\leqslant L\int^T_0 \Eof{Y^{(n)}_A(t)}\differential{t}+r_1(0)T\eqcomma                 
        \end{align*}
        for some constant $L$ since the function $r_1$ is assumed Lipschitz continuous.
        It follows from the Gr\"onwall--Bellman inequality that
        \begin{align}
        \label{eq:expectation_estimate}
            \Eof{\nY_A(t)}\leqslant r_1(0)T\myExp{LT} \text{ whence we get } \sup_{t\in[0,T]}\sup_n \Eof{\nY_A(t)}<\infty \eqstop 
        \end{align}
        Now, applying Chebyshev inequality,  one can conclude that
     \begin{align*}
    \probOf{\sup_{t\in[0,T]}\nY_A(t)\leqslant k_\epsilon}&{}\geqslant \probOf{ \tilde{Y}^{(n)}_A(T)\leqslant k_\epsilon}=1-\probOf{ \tilde{Y}^{(n)}_A(T)\geqslant k_\epsilon}\\
    &{}\geqslant 1-\frac{\int^T_0 \Eof{r_1(\nY_A(t))}\differential{t}}{k_{\epsilon}}\\
    &{}\geqslant 1-\frac{L\int^T_0 \Eof{\nY_A(t)} \differential{t}+Tr_1(0)}{k_{\epsilon}}\eqcomma             
    \end{align*}
    for some constant $L>0$. Since $\sup_{t\in[0,T]}\sup_n \Eof{\nY_A(t)}<\infty$ from \Cref{eq:expectation_estimate}, it follows from the above inequality that the process $\nY_A(t),\ t\in[0,T]$ satisfies the compact containment condition.

 Let us now define an operator $\mathcal{A}: \effectiveDomain{\mathcal{A}}\subset C_{b}(\setOfPositiveReals,\setOfReals) \mapsto C(\setOfPositiveReals\times \mathbb{N}_0, \setOfReals)$ with $\effectiveDomain{\mathcal{A}} \defeq C^{(2)}_c(\setOfPositiveReals,\setOfReals)$, the space of twice continuously differentiable functions with compact support, as follows
     \begin{align}
        \mathcal{A}f(y_1, y_2) &{}\defeq \left(r_1(y_1) - r_3(y_1)\indicator{\{0\}}{y_2}  - r_{4}(y_1) y_2  \right)\partial f(y_1) \label{eq:A_operator} 
      \eqstop  
    \end{align}        
 Now, for each  $f \in \effectiveDomain{\mathcal{A}}$, define a function $g_f$ as follows
        $
            g_f(y_1, y_2) \defeq f(y_1)\eqcomma \; \forall y_2 \in \mathbb{N}_0\eqstop 
        $
        That is, the function $g_f(y_1, y_2)$ does not depend on $y_2$. Of course, the function $g_f$ is by default in $\effectiveDomain{\mathcal{L}_n}$ for each $f \in \effectiveDomain{\mathcal{A}}$. Then, let us define the stochastic process 
        \begin{align*}
            \varepsilon_f^{(n)}(t) \defeq \int_0^t \left( \mathcal{A}f(\nY_A(s), \nY_B(s)) - \mathcal{L}_n g_f( \nY_A(s), \nY_B(s) ) \right) \differential{s} \eqstop 
        \end{align*}
        Our goal is to show that the stochastic process $\tilde{M}_f^{(n)}$ given by 
        \begin{align*}
            \tilde{M}_f^{(n)} (t) &{} = f(\nY_A(t)) - f(\nY_A(0)) - \int_{0}^{t} \mathcal{A}f(Y^{(n)}_A(s), Y^{(n)}_B(s)) \differential{s} + \varepsilon_f^{(n)}(t)
        \end{align*}
        is an $\mathcal{F}^n_t$-martingale for each $f\in \effectiveDomain{\mathcal{A}}$. It follows from the definition of the function $g_f$ that 
        \begin{align*}
            \tilde{M}_f^{(n)} (t) &{} = g_f(\nY(t)) - g_f(\nY(0)) -\int_{0}^{t}\mathcal{A}f(Y^{(n)}_A(s), Y^{(n)}_B(s)) \differential{s} + \varepsilon_f^{(n)}(t)\\ 
            &{} = g_f(\nY(t)) - g_f(\nY(0)) -\int_{0}^{t} \mathcal{L}_n g_f(\nY(s))\differential{s} 
            = M_{g_f}^{(n)}(t)\eqcomma 
        \end{align*}
        as defined in \Cref{eq:scaled_martingale}, 
        which is indeed an $\mathcal{F}^n_t$-martingale by \Cref{lem:ctmc_Dynkin_timevarying_martingale} in Appendix~\ref{sec:math_details}. Now, from \Cref{eq:A_operator}, we get 
        \begin{align*}
            (\mathcal{A}f(y_1, y_2))^2&{} \le 4\left(r_1(y_1)^2+r_3(y_1)^2+r_{4}(y_1)^4+y^4_2\right)\norm{f}_1 
            \le c_1( y^4_1+y^4_2+y^2_1+1),
        \end{align*} 
        where $\norm{f}_k=\sup_{z\in\mathbb{R}_{+},\ 0\leqslant l\leqslant k}|f^{(l)}(z)| < \infty.$ 
        Therefore, for each $f\in \effectiveDomain{\mathcal{A}}$, and for each $t\ge 0$, we have the following estimate 
        \begin{align}
            \begin{aligned}
                \Eof{ \int_{0}^{t} \absolute{\mathcal{A}f(\nY_A(s), \nY_B(s))}^2 \differential{s}  }&{}\leqslant c_1\int_{0}^{t}\left(\Eof{Y^{(n)}_A(s)^4}+\Eof{Y^{(n)}_B(s)^4}  
                \right. \\
            &{}\quad\quad\quad   \left. 
            +\Eof{Y^{(n)}_A(s)^2}\right)\differential{s} +c_1t\eqstop 
            \end{aligned}
            \label{eq:E_Af_upper_bound}
             \end{align}
        Now, note that 
   \begin{align*}
    \nY_A(t) &{}\leqslant \frac{1}{n}\int_{0}^{t} \int_{0}^{\infty}\indicator{[0, \; r^{(n)}_1(nY^{(n)}_A(u-),Y^{(n)}_B(u-))]}{v} \tilde{Q}_{A}^{(1)}(\differential{u}, \differential{v}) 
    \nonumber  \\
&{} \quad 
+ \int_{0}^{t}r^{(n)}_1(nY^{(n)}_A(u),Y^{(n)}_B(u)) \differential{u},
\end{align*}
 where $\tilde{Q}_A^{(1)}$ is the compensated \ac{PRM} corresponding to $Q^{(1)}_{A}$.   
Therefore, we get 
\begin{align*}
 \Eof{ \nY_A(t)^2} &{}\leqslant \frac{2}{n}\int_{0}^{t}  \Eof{  r_1(Y^{(n)}_A(s))} \differential{s}+2 \Eof{ \left(\int_{0}^{t} r_1(Y^{(n)}_A(s)) \differential{s}\right)^2}.
\end{align*}
It follows from \eqref{eq:expectation_estimate} and Jensen's inequality that  
\begin{align*}
 \Eof{ \nY_A(t)^2} &{}\leqslant c_2+c_3\int_{0}^{t}  \Eof{ Y^{(n)}_A(s)^2} \differential{s}.
\end{align*}
Applying the Gr\"onwall--Bellman inequality, one can conclude that
$
 \Eof{ \nY_A(t)^2} \leqslant c_2 \myExp{c_3t}.
 $
 Hence, 
 $
 \sup_n\Eof{Y^{(n)}_A(s)^2}<\infty.
 $
 Furthermore, applying the \ac{BDG} and the Jensen inequalities, one can see that 
  \begin{align*}
 \Eof{ \nY_A(t)^4} &{}\leqslant c_4\left( \Eof{\left(\frac{1}{n}\int^T_0r_1( \nY_A(s))\differential{s}\right)^2}+\Eof{\left( \int_{0}^{T}  r_1(Y^{(n)}_A(s)) \differential{s}\right)^4}\right)\\
 &{}\leqslant c_5+c_6 \int_{0}^{T} \Eof{Y^{(n)}_A(s)^4} \differential{s}. 
 \end{align*}
 Again, applying the Gr\"onwall--Bellman inequality, and taking supremum, we conclude
 $
\sup_{t\in[0,T]} \sup_n\Eof{Y^{(n)}_A(t)^4}<\infty.
 $
Now, it follows from  the \acp{SDE} in \Cref{eq:nY_sde}, the Gr\"onwall--Bellman inequality and the It\^o formula for \acp{SDE} driven \acp{PRM} (see \cite[Lemma~4.4.5]{Applebaum_2009Levy}, \cite[Theorem~5.1]{IkedaWatanabe2014Stochastic}) that
 $$
 \Eof{Y_B^{(n)}(t)}\leqslant n\lambda T \myExp{-n\mu t} \text{ and }
 \Eof{Y_B^{(n)}(t)^k}\leqslant c(k)(n\lambda T+n\mu T)\myExp{-n\mu t}\eqcomma 
  $$
  for some constants $c(k)>0$. Hence,
 $
\sup_{t\in[0,T]} \sup_n\Eof{Y^{(n)}_B(t)^4}<\infty.
 $
 Therefore, combining these estimates, we conclude  that for each $f\in \effectiveDomain{\mathcal{A}}$, and for each $t\ge 0$, 
  \begin{align}
            \sup_{n}\Eof{ \int_{0}^{t} \absolute{\mathcal{A}f(\nY_A(s), \nY_B(s))}^2 \differential{s}  }&{}< \infty.
                        \label{eq:E_Af_upper_bound_sup}
        \end{align}
Finally, we will see that
$
\lim_{n\to\infty}\Eof{\sup_{t\in[0,T]}\absolute{ \varepsilon_f^{(n)}(t)}}=0.
$ 
Indeed, since
 \begin{align}
        \begin{aligned}
 \absolute{\partial f(y_1)-n(f(y_1 +\frac{1}{n}) - f(y_1))} 
 &{}\leqslant n\int^{y_1+\frac{1}{n}}_{y_1}\absolute{f^{\prime}(y_1)-f^{\prime}(z)}\differential{z}
  \leqslant \|f\|_2\ n\int^{\frac{1}{n}}_0z\differential{z}
  \leqslant \frac{c_7}{n},
  \end{aligned} \nonumber
    \end{align}
we can estimate the expectation as follows 
  \begin{align}
        \begin{aligned}
 & \Eof{\sup_{t\in[0,T]}\absolute{ \varepsilon_f^{(n)}(t)}}\\ &{}\leqslant  \int^T_0 \Eof{ \absolute{\mathcal{A}f(\nY_A(s), \nY_B(s))-\mathcal{L}_nf(\nY_A(s), \nY_B(s)) } } \differential{s}\\
 &{}\leqslant \int^T_0\Eof{r_1(\nY_A(s))\absolute{\partial f(\nY_A(s))-n(f(\nY_A(s) +\frac{1}{n}) - f(\nY_A(s)))} }\differential{s}\\
 &{}\quad + \int^T_0\Eof{r_3(\nY_A(s))\absolute{- \partial f(\nY_A(s))-n(f(\nY_A(s) -\frac{1}{n}) - f(\nY_A(s)))} }\differential{s}\\ 
 &{}\quad + \int^T_0\Eof{r_{4}(\nY_A(s))\nY_B(s)\absolute{- \partial f(\nY_A(s))-n(f(\nY_A(s)-\frac{1}{n}) - f(\nY_A(s)))} }\differential{s}\\
 &{}\leqslant \frac{c_7}{n}\int^T_0\Eof{r_1(\nY_A(s))}\differential{s}+\frac{c_7}{n}\int^T_0\Eof{r_3(\nY_A(s))}\differential{s}
 +\frac{c_7}{n}\int^T_0\Eof{r_{4}(\nY_A(s))\nY_B(s)}\differential{s}
 \\
 &{}
 \leqslant \frac{c_8}{n}\eqcomma \nonumber 
   \end{aligned}
    \end{align}
for some constant $c_8>0$ since, by the Cauchy--Schwarz inequality, we have 
$$
\Eof{r_{4}(\nY_A(s))\nY_B(s)}\leqslant\left(\Eof{r_{4}(\nY_A(s))^2}\right)^{\frac{1}{2}}\left(\Eof{\nY_B(s)^2}\right)^{\frac{1}{2}}\eqstop 
$$
Hence, we obtain that 
 \begin{align}
 \lim_{n\to\infty}\Eof{\sup_{t\in[0,T]}\absolute{ \varepsilon_f^{(n)}(t)}}=0.
 \label{eq:zero_lim}
 \end{align} 
 The compact containment condition for $\nY_A$, \eqref{eq:E_Af_upper_bound_sup} and  \eqref{eq:zero_lim} imply the relative compactness of $\nY_A$ in $D([0, T], \setOfPositiveReals)$ \cite[Chapter 3]{Ethier:1986:MPC}. 
 
 Let us now check the relative compactness of the sequence of  random measures $\{\Gamma_n : n\ge 1\}$. The most straightforward way to verify the relative compactness of a collection of random measures is to verify the relative compactness of the corresponding collection of mean measures (also called intensities).  To this end, we will show that for each 
 $\delta>0$, there exists a constant $M_0>0$ such that for any $M_0\leqslant M<\infty$
 $$
 \inf_n \Eof{\Gamma_n([0,t]\times [0,M])}\geqslant t(1-\delta).
 $$
 Indeed, it follows in a straightforward manner that 
\begin{align}
        \begin{aligned}
            \Eof{\Gamma_n([0,t]\times [0,M])}&{}=\Eof{\int^t_0\indicator{[0,M]}{\nY_B}(s)\differential{s}}
            = \int^t_0\probOf{Y_B^{(n)}(s)\leqslant M} \differential{s}\\
            &{} \geqslant t-\int^t_0\probOf{Y_B^{(n)}(s)\geqslant M+1} \differential{s}
            \\
            &{}
            \geqslant t-\int^t_0\frac{ \Eof{(Y_B^{(n)}(s))}}{M+1}\differential{s} 
            \geqslant t\left(1-\frac{c}{M+1}\right)
        \end{aligned}
    \end{align}
    for some constant $c>0$ in the light of our previous estimates (as well as \Cref{lem:MGF_hitting_times}). 
    Hence, the sequence of random measures $\{\Gamma_n\}_{n \ge 1}$ is relatively compact in $\mathcal{M}_1([0, T]\times \mathbb{N}_{0})$ by virtue of \cite[Lemma 1.3]{Kurtz1992Averaging}. 
    
    Joint relative compactness of the vector $(\nY_A, \Gamma_n)$ in the space $D([0, T], \setOfPositiveReals)\times \mathcal{M}_1([0, T]\times \mathbb{N}_{0})$  follows from the relative compactness of the individual components \cite[Problems 5.9, p.65]{Billingsley1999Convergence}. The topologies on the respective spaces $D([0, T], \setOfPositiveReals)$ and $ \mathcal{M}_1([0, T]\times \mathbb{N}_{0})$ are mentioned in notational conventions in \Cref{sec:intro}. 

    \paragraph{Identification of the limit} 
    Since the sequence $\{(\nY_A, \Gamma_n) : n \ge 1\}$ is relatively compact, we can find a convergent subsequence. Let $(y_A,\Gamma)$ be a limit point of $(Y_A^{(n)},\Gamma_n)$. 
    Define the filtration 
    $
    \mathcal{H}_t=\sigma\{Y_A(s),\ \Gamma([0,s]\times H),\ s\leqslant t,\ H\in\mathcal{B}(\mathbb{N}_0)\}.
    $
    Then, the fact that
     \begin{align*}
            f(y_A(t)) - f(y_A(0)) - \int_{0}^{t} \int_{\setOfNaturals_0} \mathcal{A} f(y_A(s), z) \Gamma(\differential{s}\times \differential{z})
        \end{align*}
        is an $\mathcal{H}_t$-martingale, for any limit point $(y_A,\Gamma)$ of $(Y_A^{(n)},\Gamma_n)$, and all $f \in \effectiveDomain{\mathcal{A}}$ follows directly from \cite[Theorem 2.1]{Kurtz1992Averaging}, which we recollect as \Cref{thm:stochastic_averaging} in Appendix~\ref{sec:math_details} for the sake of completeness and the readers' convenience. 

        Now, in order to characterise the limit points,  define the operator $\mathcal{B} : \effectiveDomain{\mathcal{B}} \defeq  C_b(\setOfNaturals_0,\mathbb{R})\mapsto C_b(\setOfPositiveReals\times \setOfNaturals_0,\mathbb{R})$ as follows 
        $
            \mathcal{B} p(y_1, y_2) \defeq \mathcal{B}_{y_1} p(y_2) \eqcomma 
            $
        where the family of operators $\mathcal{B}_{y_1}$ are defined in \Cref{eq:frozen_generator} for all $y_1 \in \setOfPositiveReals$. Also, for $p \in \effectiveDomain{\mathcal{B}}$, define $h_p $ as $h_p(y_1, y_2) = p(y_2)$ for all $y_2\in\mathbb{R}_+$. That is, $h_p(y_1, y_2)$ does not depend on $y_1$. Now, define the stochastic process 
        \begin{align*}
            \delta_{p}^{(n)}(t) \defeq \int_{0}^{t} \left( n \mathcal{B} p(\nY_A(u), \nY_B(u)) - \mathcal{L}_n h_p( \nY_A(u), \nY_B(u))\right)\differential{u}\eqcomma 
        \end{align*}
        for  $p \in \effectiveDomain{\mathcal{B}}$. Our goal is to show that the stochastic process 
        \begin{align*}
            \hat{M}_p^{(n)}(t) \defeq p(\nY_B(t)) - p(\nY_B(0)) -\int_{0}^{t} n  \mathcal{B} p(\nY_A(u), \nY_B(u))\differential{u} + \delta_{p}^{(n)}(t)\eqcomma \quad t\ge 0\eqcomma 
        \end{align*}
        is an $\history{t}^n$-martingale for all $p \in \effectiveDomain{\mathcal{B}}$. Indeed, 
        \begin{align*}
            \hat{M}_p^{(n)} (t) &{} = h_p(\nY(t)) - h_p(\nY(0)) -\int_{0}^{t}n \mathcal{B}p(nY_A(s), nY_B(s)) \differential{s} + \delta_p^{(n)}(t)\\ 
            &{} = h_g(\nY(t)) - h_g(\nY(0)) -\int_{0}^{t} \mathcal{L}_n h_p(\nY(s))\differential{s} 
            = M_{h_p}^{(n)}(t)\eqcomma 
        \end{align*}
        from \Cref{eq:scaled_martingale}, which is indeed an $\history{t}^n$-martingale by \Cref{lem:ctmc_Dynkin_timevarying_martingale} in Appendix~\ref{sec:math_details}. Moreover, it is straightforward to verify that 
        \begin{align*}
            \lim_{n\to \infty} \Eof{\sup_{s \le t}  n^{-1} \absolute{ \delta_{p}^{(n)}(s) } } = 0\eqcomma \text{ for all } t\ge 0\eqstop 
        \end{align*}
        Since the Markov process corresponding to the frozen generator $\mathcal{B}_{y_1}$ admits a unique stationary distribution $\pi$ defined in \Cref{eq:stationary_distribution}, it follows from \cite[Example 2.3]{Kurtz1992Averaging} that the limiting measure $\Gamma$ on $\setOfPositiveReals\times \setOfNaturals_0$ can be represented as 
          \begin{align*}
            \int_{0}^{t} \int_{\setOfNaturals_{0}} f(s, z) \Gamma(\differential{s}\times \differential{z}) = \int_{0}^{t} \int_{\mathbb{N}_0}f(s, z)\pi(\differential{z})\times\differential{s}= \int_{0}^{t}\left( \sum_{z=0}^{\infty } f(s, z)\pi(z)\right)\differential{s} \eqstop 
        \end{align*}
        Hence, the martingale associated with the limit point $(y_A, \Gamma)$ is given by
         \begin{align*}
            \mathcal{M}_f(t)\defeq  f(y_A(t)) - f(y_A(0)) - \int_{0}^{t} \int_{\setOfNaturals_0} \mathcal{A} f(y_A(s), z) \Gamma(\differential{s}\times \differential{z})\\
            =f(y_A(t)) - f(y_A(0)) - \int_{0}^{t} \int_{\mathbb{N}_0}\mathcal{A}f(y_A(s), z)\pi(\differential{z})\times\differential{s}.
        \end{align*}
        Note that  $\mathcal{M}_f$ is a square integrable martingale with $\langle\mathcal{M}_f\rangle (t)=0$ for all $t\ge 0$.  It implies that
        $\mathcal{M}_f (t)=0$ almost surely for all $t\ge 0$.  Hence, we get 
        \begin{align*}
            f(y_A(t)) - f(y_A(0)) &{} =\int_{0}^{t} \int_{\mathbb{N}_0}\mathcal{A}f(y_A(s), z)\pi(\differential{z}) \differential{s}\\
        &{} = \int_{0}^{t} \left(\sum_{z=0}^{\infty } \left(r_1(y_A(s))-r_3(y_A(s))\indicator{z}{0}-r_{4}(y_A(s))z\right)\pi(z)\right) \\
        &{}\quad \quad \times \partial f(y_A(s)) \differential{s}\\
        &{} = \int_0^t \left(r_1(y_A(s))-r_3(y_A(s)) e^{-m(y_A(s))} -r_4(y_A(s))m(y_A(s)) \right) \\
        &{}\quad \quad \times \partial f(y_A(s)) \differential{s}\eqstop 
        \end{align*}
        In fact, the limit $y_A$ lies in $C([0, T], \setOfPositiveReals)$ almost surely by virtue of \cite[Theorem 13.4]{Billingsley1999Convergence}. 
        Therefore, the weak convergence of the sequence  $\{ \nY_A : n\ge 1\}$ to $y_A$ follows from  the fundamental theorem of calculus, and the existence and the uniqueness of the solutions of the \ac{ODE} in \Cref{eq:limit_ODE}. Since the limit $y_A$ is deterministic, the convergence also holds in probability. 
%
   \end{proof}

    \begin{figure}[t!]
         \begin{subfigure}[b]{1\columnwidth}
        \includegraphics[width=0.49\columnwidth,trim={0cm 0cm 1cm 1cm},clip ]{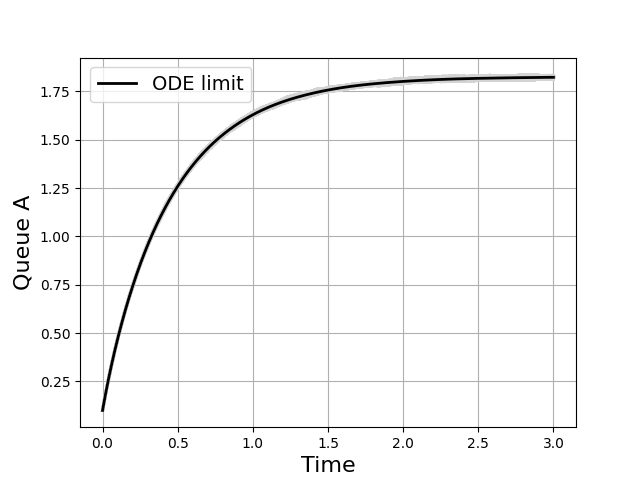}
        \includegraphics[width=0.49\columnwidth,trim={0cm 0cm 1cm 1cm},clip]{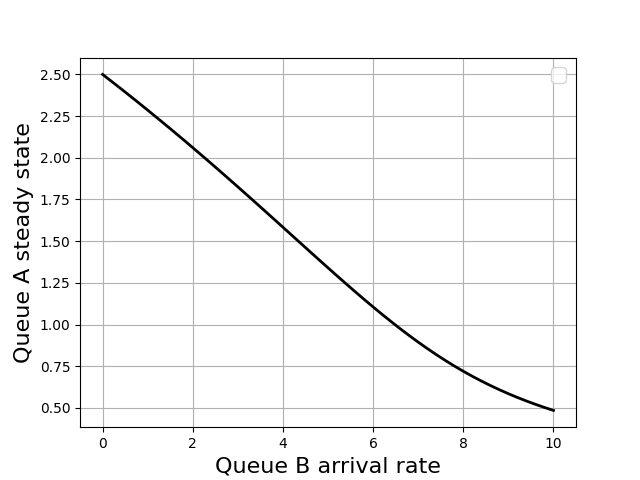}
        \vspace{-10pt}
        \caption{}
        \label{fig: speedup}
        \end{subfigure}
            \begin{subfigure}[b]{1\columnwidth}
        \includegraphics[width=0.49\columnwidth,trim={0cm 0cm 1cm 1cm},clip]{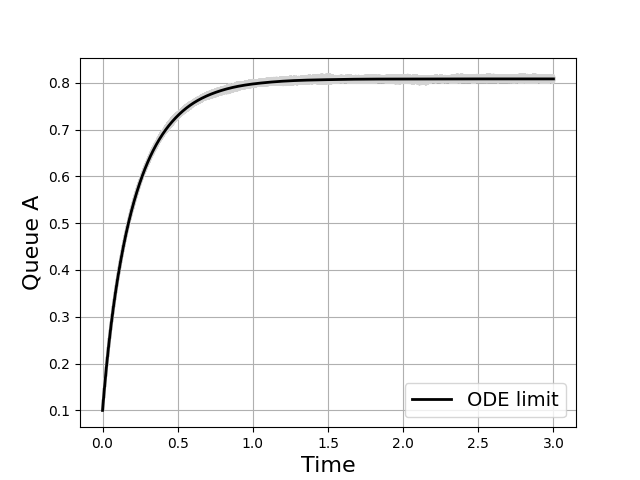}
        \includegraphics[width=0.49\columnwidth,trim={0cm 0cm 1cm 1cm},clip]{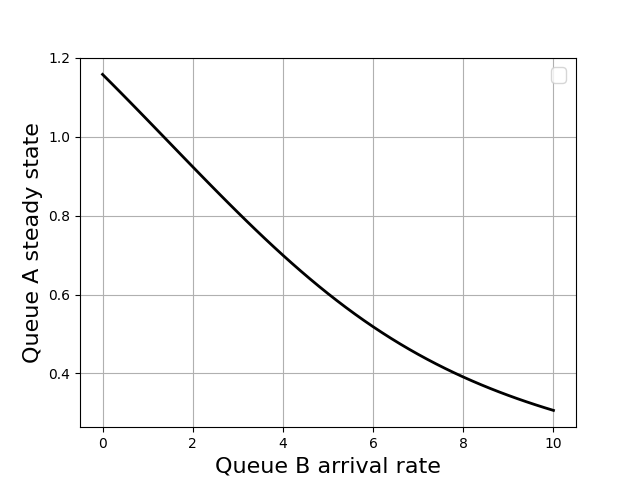}
        \vspace{-10pt}
        \caption{}
        \label{fig: regulated_arrivals}
        \end{subfigure}
        \vspace{-25pt}
        \caption{\textbf{(Left)}  Accuracy of the \ac{FLLN} approximation for the scaled queue lengths $\nY_A$ by the deterministic function $y_A$ solving the \ac{ODE} in \Cref{eq:limit_ODE} in \Cref{thm:relative_compactness}. The stochastic simulations are performed using the standard Doob--Gillespie's algorithm. A total of $100$ trajectories are shown in this figure. \textbf{(Right)} Steady state  values of the scaled queue length at queue~$A$ as a function of the arrival rate $\lambda_B$ at the queue $B$. The queue A dependent propensity functions are: (a) $r_1=\lambda_A, r_3(y_1)=\mu_A y_1, r_4(y_1)=M \mu_A y_1$ (b) $r_1=\lambda_A/(y_1+1), r_3(y_1)=\mu_A y_1, r_4(y_1)=M \mu_A y_1 $. The parameter values used: $n=10^5, y_A(0)=0, \lambda_A = 5.0, \lambda_B=3.0, \mu_A=2.0, \mu_B=2.0$, and the speedup factor $M=2.0$. 
        } 
        \label{fig:FLLN}
    \end{figure}
   

    \Cref{fig:FLLN} shows that the \ac{FLLN} is accurate in that exactly simulated trajectories of $\nY_A$ with a finite but large $n$ are close to the deterministic limit $y_A$ solving the \ac{ODE} \eqref{eq:limit_ODE}. The right side panel of \Cref{fig:FLLN} is interesting from a communication engineering perspective and shows that the steady state of the queue $A$ can be controlled as a function of the arrival rate $\lambda_B$ at the queue $B$.

    \section{\acl{FCLT}}
    \label{sec:FCLT}
    We now turn our attention to the \ac{FCLT}. We will show that the fluctuations of the scaled process $\nY_A$ around the deterministic trajectory $y_A$ are Gaussian. To this end, let us  define the scaled stochastic process $W_n$ as follows
   $$
   W_n(t)=\sqrt{n}(\nY_A(t)-y_A(t)),\ t\in[0,T].
   $$
   In order to prove an \ac{FCLT} for the sequence of  stochastic processes $W_n$, we first need to find the solution $F(y_1, \cdot)$ to the Poisson equation  (see \cite{Glynn1996Liapounov,Costa2003Poisson,Kang:2014:CLT}) 
   \begin{align}
        \mathcal{B}_{y_1}F(y_1, \cdot)(y_2) = - h_{y_1}(y_2)
        \label{eq:Poisson_eqn}
   \end{align}
   for each $y_1 \in \setOfPositiveReals$, where the frozen generator $\mathcal{B}_{y_1}$ is defined in \Cref{eq:frozen_generator}, and the function $h_{y_1}(\cdot) : \setOfNonnegativeIntegers \mapsto \setOfReals$ is given by
   \begin{align*}
    h_{y_1}(y_2)\defeq r_3(y_1)(\indicator{\{0\}}{y_2}- \myExp{-m(y_1)})+r_{4}(y_1)(y_2-m(y_1)) \eqstop 
   \end{align*}
   The following lemma provides the solution to the Poisson equation.
   \begin{myLemma}
    The solution $F$ to the Poisson equation \eqref{eq:Poisson_eqn} is given by 
    \begin{align*}
        F(y_1, y_2) = u_1(y_1)\indicator{\{0\}}{y_2}+u_2(y_1)y_2\eqcomma
    \end{align*}
    where the functions $u_1$ and $u_2$ are defined as
    \begin{align*}
        u_1(y_1) &{} = \frac{1}{\lambda}\left[\lambda u_2(y_1) + r_3(y_1)(1- \myExp{m(1)}) - r_4(y_1)m(y_1)\right]\eqcomma \\
        u_2(y_1) &{} = \frac{1}{ \lambda^2}\left[\lambda r_3(y_1) \myExp{-m(y_1)} -\lambda r_4(y_1)(1- m(y_1))  \right. \\
        &{}\quad\quad \quad  \left.
            +r_3(y_1)(1- \myExp{-m(y_1)}) -r_4(y_1) m(y_1) \right] \eqstop 
    \end{align*}
    \label{lem:Poisson_eqn} 
   \end{myLemma}
   \begin{proof}[Proof of \Cref{lem:Poisson_eqn}]
        Assume that the solution has the form  $$F(y_1, y_2) = u_1(y_1)\indicator{\{0\}}{y_2}+u_2(y_1)y_2\eqstop $$
        It is enough to solve for  
        $
            F(y_1, 0)  = u_1(y_1)\eqcomma $ and $  F(y_1, 1) = u_2(y_1)\eqstop 
        $
        Then, equating coefficients on both sides of the Poisson equation \eqref{eq:Poisson_eqn}, we get the following two linear equations in $u_1$ and $u_2$:
        \begin{align*} 
            \lambda(u_2(y_1)-u_1(y_1)) &{} =-r_3(y_1)(1- \myExp{-m(y_1)})
            +r_{4}(y_1)m(y_1)\eqcomma \\
            \lambda u_2(y_1)+(r_{4}(y_1)+\mu)(u_1(y_1)-u_2(y_1))&{} =r_3(y_1)\myExp{-m(y_1)} 
            -r_{4}(y_1)(1-m(y_1))\eqstop 
        \end{align*} 
        Solving the above equations, we get the expressions for $u_1$ and $u_2$ 
        as stated in the lemma. It is straightforward to verify that the solution $F(y_1, \cdot)$ thus obtained satisfies the Poisson equation \eqref{eq:Poisson_eqn} for all $y_2\in \setOfNonnegativeIntegers$. 
   \end{proof}

   We are now ready to state and prove the \ac{FCLT} for the scaled process $W_n$, which is the main result of this section. 

       \begin{myTheorem}[\acl{FCLT}]
        Assume $W_n (0) = 0$ for all $n\ge 1$, and the conditions of \Cref{thm:relative_compactness}. Further, assume that the functions $r_1, r_3$, and $r_4$ are twice continuously differentiable. Then, the scaled stochastic  process $W_n$ converges weakly to a continuous Gaussian semimartingale $W$ 
        as $n\to\infty$ where the process $W$ lies in $C([0,T], \setOfReals),$ and satisfies the following \ac{SDE} (written in the integral form)
        \begin{align}
            W(t) = \int_{0}^{t} \partial G(y_A(s))W(s) \differential{s} + \int_{0}^{t} \sqrt{ \sigma_F(s) }\differential{B(s)}\eqcomma \label{eq:clt_sde}
        \end{align}
        where the process  $B$ is a standard Brownian motion, the function $G$ is given by \Cref{eq:limit_ODE}, the positive  function $\sigma_F$ is given by
        \begin{align*}
            \sigma_F(s) &{} \defeq \left(r_1(y_A(s))+r_3(y_A(s))\myExp{-m(y_A(s))}\right)\\ 
            &{}\quad +  \left[ (u_1(y_A(s))-u_2(y_A(s)) + 1)^2 \pi(1) 
             + (1- u_2(y_A(s)))^2 (m(y_A(s))-\pi(1)) \right]r_{4}(y_A(s))\\
            &{}\quad + \lambda \left[(u_2(y_A(s)) - u_1(y_A(s))   )^2 \pi(0) + u_2(y_A(s))^2 \pi([1, \infty) ) \right] \\
            &{}\quad 
            + \mu\left[ (u_1(y_A(s))-u_2(y_A(s)))^2 \pi(1) + u_2(y_A(s))^2 (m(y_A(s))-\pi(1)) \right] \eqcomma
        \end{align*}
        the function $y_A$ solves the \ac{ODE} in \Cref{eq:limit_ODE},
        and the functions $u_1, u_2$ are given in the definition of $F(y_1, \cdot)$,  the solution to the Poisson equation \eqref{eq:Poisson_eqn} as given in \Cref{lem:Poisson_eqn}. 
       \label{thm:clt}
    \end{myTheorem}
     \begin{proof}[Proof of \cref{thm:clt}]
     Note that the function $y_A$ satisfies 
     $$
     y_A(t)=\int^t_0r_1(y_A(s))\differential{s}-\int^t_0r_3(y_A(s))e^{-m(y_A(s))}\differential{s}-\int^t_0r_{4}(y_A(s))m(y_A(s))\differential{s}.
     $$
     Therefore, from the \ac{SDE} representation of $\nY_A$ in \Cref{eq:nY_sde}, we have 
        \begin{align*}
        W^{(n)}_A(t)&{}=M^{(n)}_A(t)+\sqrt{n}\int^t_0(r_1(\nY_A(s))-r_1(y_A(s)))\differential{s}\\
        &{}\quad -\sqrt{n}\int^t_0(r_3(\nY_A(s))\indicator{\{0\}}{\nY_B(s)}-r_3(y_A(s))\myExp{-m(y_A(s))})\differential{s}\\
          &{}\quad -\sqrt{n}\int^t_0(r_{4}(\nY_A(s))\nY_B(s)-r_{4}(y_A(s))m(y_A(s))\differential{s}\\
          &{}=M^{(n)}_A(t)+\sqrt{n}\int^t_0(r_1(\nY_A(s))-r_1(y_A(s)))\differential{s}\\
          &{}\quad -\sqrt{n}\int^t_0(r_3(\nY_A(s))e^{-m(Y^{(n)}_A(s))}-r_3(y_A(s))e^{-m(y_A(s))})\differential{s}\\
           &{}\quad -\sqrt{n}\int^t_0(r_{4}(\nY_A(s))m(Y^{(n)}_A(s))-r_{4}(y_A(s))m(y_A(s))\differential{s}\\         
          &{}\quad -\sqrt{n}\int^t_0(r_3(\nY_A(s))))\indicator{\{0\}}{\nY_B(s)}-r_3(\nY_A(s))e^{-m(Y^{(n)}_A(s))})\differential{s}\\
          &{}\quad -\sqrt{n}\int^t_0(r_{4}(Y^{(n)}_A(s))Y^{(n)}_B(s)-r_{4}(\nY_A(s))m(Y^{(n)}_A(s)))\differential{s}\\
          &{} =  M^{(n)}_A(t)+\sqrt{n}\int^t_0(G(Y^{(n)}_A)-G(y_A(s)))\differential{s}-\sqrt{n}\int^t_0h_{Y^{(n)}_A(s)}(Y^{(n)}_B(s))\differential{s}\eqcomma 
          \end{align*} 
        where the function $h$ is given by
        \begin{align*}
            h_{y_1}(y_2)=r_3(y_1)(\indicator{\{0\}}{y_2}- \myExp{-m(y_1)})+r_{4}(y_1)(y_2-m(y_1)),
        \end{align*}
        the stochastic process $M^{(n)}_A$ is defined as
        \begin{align*}
            M^{(n)}_A(t) &{} = \frac{1}{\sqrt{n}}\int_{0}^{t} \int_{0}^{\infty}\indicator{[0, \; r^{(n)}_1(nY^{(n)}_A(u-),Y^{(n)}_B(u-))]}{v} \tilde{Q}_{A}^{(1)}(\differential{u}, \differential{v}) \nonumber  \\
            &{} \quad -\frac{1}{\sqrt{n}} \int_{0}^{t} \int_{0}^{\infty}\indicator{[0, \; r^{(n)}_3(nY^{(n)}_A(u-),Y^{(n)}_B(u-))]}{v} \tilde{Q}_A^{(2)}(\differential{u}, \differential{v}) \nonumber   \\
            &{} \quad - \frac{1}{\sqrt{n}}\int_0^t \int_0^\infty \indicator{[0, \; r^{(n)}_4(nY^{(n)}_A(u-),Y^{(n)}_B(u-))]}{v}\tilde{Q}_{AB}^{(1)}(\differential{u}, \differential{v})\eqcomma 
        \end{align*}
        and $    \tilde{Q}_{A}^{(1)}, \tilde{Q}_{A}^{(2)}$, and $\tilde{Q}_{AB}^{(1)}$ are the compensated \acp{PRM} associated with $Q_A^{(1)}, Q_A^{(2)}$, and $Q_{AB}^{(1)}$ respectively.  

On the other hand,  the generator of the process $Y^{(n)}$ can be represented as follows
$$
\mathcal{L}_nf(y_1,y_2)=\delta^{(n)}_f(y_1,y_2)+n\mathcal{B}_{y_1}f(y_1,\cdot)(y_2),
$$
where the function $\delta^{(n)}_f$ is given by 
\begin{align}
 \begin{aligned}
 \delta^{(n)}_f(y_1,y_2) &{}\defeq nr_1(y_1)\left(f(y_1+\frac{1}{n},y_2)-f(y_1,y_2)\right)\\
 &{}\quad +nr_3(y_1)\indicator{\{0\}}{y_2}\left(f(y_1-\frac{1}{n},y_2)-f(y_1,y_2)\right)\\
 &{}\quad + nr_4(y_1)y_2\left(f(y_1-\frac{1}{n},y_2-1)-f(y_1,y_2-1)\right)
       \end{aligned}
       \label{eq:delta_f_n}
\end{align}
 and the generator $\mathcal{B}_{y_1}$ is defined in  \Cref{eq:frozen_generator}. It is immediate  that for each $y_1\in\setOfPositiveReals$ and $y_2\in\setOfNonnegativeIntegers$, we have 
  \begin{align*}
   \delta^{(n)}_f(y_1,y_2)\to r_1(y_1)\frac{\partial f(y_1,y_2)}{\partial y_1}
 -r_3(y_1)\indicator{0}{y_2}\frac{\partial f(y_1,y_2)}{\partial y_1}
 -r_4(y_1)y_2 \frac{\partial f(y_1,y_2-1)}{\partial y_1}
  \end{align*}
as $n\to\infty.$ Now, let $F(y_1, \cdot)$ be the solution of the Poisson equation 
$
    \mathcal{B}_{y_1}F(y_1,y_2)=-h_{y_1}(y_2)
$
from \Cref{lem:Poisson_eqn}. Let us apply It\^o's lemma (see \cite[Lemma~4.4.5]{Applebaum_2009Levy}, \cite[Theorem~5.1]{IkedaWatanabe2014Stochastic}) to the \acp{SDE} in \Cref{eq:nY_sde} with the function $F$ to get
  \begin{align*}
  F(Y^{(n)}(t))-F(Y^{(n)}(0))&{}=\int^t_0\mathcal{L}_n F(Y^{(n)}(s))\differential{s}+M^{(n)}_F(t)\\  
  &{}=\int^t_0 n \mathcal{B}_{Y^{(n)}_A(s)}F(Y^{(n)}_A(s),\cdot)(Y^{(n)}_B(s))\differential{s}
  \\
  &{}\quad 
  +\int^t_0\delta^{(n)}_F(Y^{(n)}(s))\differential{s}+M^{(n)}_F(t)\\
   &{}=-n \int^t_0h_{Y^{(n)}_A(s)}(Y^{(n)}_B(s))\differential{s} 
  +\int^t_0\delta^{(n)}_F(Y^{(n)}(s))\differential{s}+M^{(n)}_F(t),
        \end{align*}
 where the martingale $M^{(n)}_F(t)$ is given by
 \begin{align*}
        \begin{aligned}
             M^{(n)}_F(t)&{} = \int_{0}^{t} \int_{0}^{\infty}\left(F(Y^{(n)}_A(s-)+\frac{1}{n},Y^{(n)}_B(s-))-F(Y^{(n)}(s-))\right) \\
             &{}\quad\quad \quad  
             \times \indicator{[0, \; n r_1(Y^{(n)}_A(s-)) ]}{v} \tilde{Q}_{A}^{(1)}(\differential{s}, \differential{v}) \\
            &{} \quad +\int_{0}^{t} \int_{0}^{\infty}\left(F(Y^{(n)}_A(s-)-\frac{1}{n},Y^{(n)}_B(s-))-F(Y^{(n)}(s-))\right) \\
            &{}\quad \quad \quad 
            \times \indicator{[0, \; n r_3(Y^{(n)}_A(s-))\indicator{\{0\}}{Y^{(n)}_B(s-)}]}{v} \tilde{Q}_{A}^{(2)}(\differential{s}, \differential{v})  \\
            &{} \quad +\int_{0}^{t} \int_{0}^{\infty}\left(F(Y^{(n)}_A(s-)-\frac{1}{n},Y^{(n)}_B(s-)-1)-F(Y^{(n)}(s-))\right)\\
            &{}\quad \quad \quad 
            \times \indicator{[0, \; n r_4(Y^{(n)}_A(s-)) Y^{(n)}_B(s-)]}{v} \tilde{Q}_{AB}^{(1)}(\differential{s}, \differential{v})  \\
            &{} \quad +\int_{0}^{t} \int_{0}^{\infty}\left(F(Y^{(n)}_A(s-),Y^{(n)}_B(s-)+1)-F(Y^{(n)}(s-))\right)
            \indicator{[0, \; n\lambda]}{v} \tilde{Q}_{B}^{(1)}(\differential{s}, \differential{v})  
            \\
            &{} \quad 
            +\int_{0}^{t} \int_{0}^{\infty}\left(F(Y^{(n)}_A(s-),Y^{(n)}_B(s-)-1)-F(Y^{(n)}(s-))\right)
            \\
            &{}\quad \quad \quad 
            \times 
            \indicator{[0, \; n\mu Y^{(n)}_B(s-)]}{v} \tilde{Q}_{B}^{(2)}(\differential{s}, \differential{v})\eqcomma 
        \end{aligned}
    \end{align*}
    where $\tilde{Q}_{A}^{(1)}, \tilde{Q}_{A}^{(2)}, \tilde{Q}_{AB}^{(1)}, \tilde{Q}_{B}^{(1)}$, and $\tilde{Q}_{B}^{(2)}$ are the compensated \acp{PRM} associated with $Q_A^{(1)}, Q_A^{(2)}, Q_{AB}^{(1)}, Q_B^{(1)}$, and $Q_B^{(2)}$ respectively. Hence, we get 
   \begin{align*}
  \sqrt{n}\int^t_0h_{Y^{(n)}_A(s)}(Y^{(n)}_B(s))\differential{s} &{}=\frac{1}{\sqrt{n}}(F(Y^{(n)}(t))-F(Y^{(n)}(0)))\\
   &{}\quad -\frac{1}{\sqrt{n}}\int^t_0\delta^{(n)}_F(Y^{(n)}(s))\differential{s}-\frac{1}{\sqrt{n}}M^{(n)}_F(t).
    \end{align*}  
Putting the above expression in the representation of $W^{(n)}_A$, and using the Taylor's theorem with integral remainder, we get
 \begin{align} 
 W_n(t)
 &{}=M^{(n)}_A(t)+\sqrt{n}\int^t_0(G(Y^{(n)}_A)-G(y_A(s)))\differential{s}-\frac{1}{\sqrt{n}}(F(Y^{(n)}(t))-F(Y^{(n)}(0))) \nonumber \\
   &{}\quad +\frac{1}{\sqrt{n}}\int^t_0\delta^{(n)}_F(Y^{(n)}(s))\differential{s}+\frac{1}{\sqrt{n}}M^{(n)}_F(t) \nonumber \\
  &{}=\tilde{M}^{(n)}(t)+\sqrt{n}\int^t_0(G(Y^{(n)}_A)-G(y_A(s)))\differential{s}-\frac{1}{\sqrt{n}}(F(Y^{(n)}(t))-F(Y^{(n)}(0))) \nonumber \\ 
  &{}\quad +\frac{1}{\sqrt{n}}\int^t_0\delta^{(n)}_F(Y^{(n)}(s))\differential{s} \nonumber \\
  &{} = \tilde{M}^{(n)}(t)+\int^t_0\left(\partial G(y_A(s))W_n(s)+ \frac{\sqrt{n}}{2}\int^{Y^{(n)}_A(s)}_{y_A(s)}\partial^2 G(r)(Y^{(n)}_A(s)-r)\differential{r}\right)\differential{s} \nonumber \\
  &{}\quad 
  -\frac{1}{\sqrt{n}}(F(Y^{(n)}(t))-F(Y^{(n)}(0)))+\frac{1}{\sqrt{n}}\int^t_0\delta^{(n)}_F(Y^{(n)}(s))\differential{s}  \label{eq:main_repres}
  \end{align}
  where $\tilde{M}^{(n)}(t)=M^{(n)}_A(t)+\frac{1}{\sqrt{n}}M^{(n)}_F(t)$ is a square-integrable martingale. 
  We will study each term on the right-hand side of \Cref{eq:main_repres} separately. 

     
\paragraph{Convergence of the martingale term $\tilde{M}^{(n)}$}
First, note that the expected value of the maximum jump in $\tilde{M}^{(n)}$ is asymptotically negligible, \ie, 
 \begin{align*}
  \Eof{\sup_{0<t\leqslant T}\absolute{\tilde{M}^{(n)}(t)-\tilde{M}^{(n)}(t-)}}\leqslant\frac{1}{\sqrt{n}}\to0,\ \text{ as }n\to\infty.
   \end{align*} 
  Moreover, the quadratic variation of $\tilde{M}^{(n)}$ has the following form
   \begin{align*}
    \predictableVariation{\tilde{M}^{(n)}}(t) &{}= \int_{0}^{t}\left(F\left(Y^{(n)}_A(s)+\frac{1}{n},Y^{(n)}_B(s)\right)-F(Y^{(n)}(s))+1\right)^2r_1(\nY_A(s))\differential{s} \\
            &{}\quad +\int_{0}^{t}\left(F\left(Y^{(n)}_A(s)-\frac{1}{n},Y^{(n)}_B(s)\right)-F(Y^{(n)}(s))+1\right)^2r_3(\nY_A(s))\\
            &{}\quad \quad \quad 
            \times \indicator{\{0\}}{\nY_B(s)} \differential{s} \\
            &{}\quad +\int_{0}^{t}\left(F\left(Y^{(n)}_A(s)-\frac{1}{n},Y^{(n)}_B(s)-1\right)-F(Y^{(n)}(s))+1\right)^2 \\
            &{}\quad \quad \quad 
            \times r_{4}(\nY_A(s))\nY_B(s)\differential{s} \\
            &{}\quad +\int_{0}^{t}\lambda \left(F\left(Y^{(n)}_A(s),Y^{(n)}_B(s)+1\right)-F(Y^{(n)}(s))\right)^2 \differential{s}\\
            &{}\quad +\int_{0}^{t}\mu \left(F\left(Y^{(n)}_A(s),Y^{(n)}_B(s)-1\right)-F(Y^{(n)}(s))\right)^2 Y^{(n)}_B(s)\differential{s}\\
            &{}= \int_{0}^{t}\int_{\mathbb{N}_0}\left(F\left(Y^{(n)}_A(s)+\frac{1}{n},z\right)-F(Y^{(n)}_A(s),z)+1\right)^2 \\
            &{}\quad \quad \quad 
            \times r_1(\nY_A(s))\Gamma_n(\differential{s}\times \differential{z})\\
            &{}\quad +\int_{0}^{t}\int_{\mathbb{N}_0}\left(F\left(Y^{(n)}_A(s)-\frac{1}{n},z\right)-F(Y^{(n)}_A(s),z)+1\right)^2\\
            &{}\quad \quad \quad 
            \times r_3(\nY_A(s))\indicator{\{0\}}{z}\Gamma_n(\differential{s}\times \differential{z}) \\
            &{}\quad +\int_{0}^{t}\int_{\mathbb{N}_0}\left(F\left(Y^{(n)}_A(s)-\frac{1}{n},z-1\right)-F(Y^{(n)}_A(s),z)+1\right)^2\\
            &{}\quad \quad \quad 
            \times r_{4}(\nY_A(s))z\Gamma_n(\differential{s}\times \differential{z}) \\
            &{}\quad +\int_{0}^{t}\int_{\mathbb{N}_0}\lambda \left(F\left(Y^{(n)}_A(s),z+1\right)-F(Y^{(n)}_A(s),z))\right)^2 \Gamma_n(\differential{s}\times \differential{z})\\
            &{}\quad +\int_{0}^{t}\int_{\mathbb{N}_0}\mu \left(F\left(Y^{(n)}_A(s),z-1\right)-F(Y^{(n)}_A(s),z))\right)^2 z\Gamma_n(\differential{s}\times \differential{z})\eqstop                       
    \end{align*}  
    By virtue of  \Cref{thm:relative_compactness}, we can conclude that 
    \begin{align*}
        \predictableVariation{\tilde{M}^{(n)}}(t) \ConvInProb \Sigma_F(t) \eqcomma 
    \end{align*}
    where $\Sigma_F(t)$ is given by 
    \begin{align*}
        \Sigma_F(t) &{} = \int^t_0\left(r_1(y_A(s))+r_3(y_A(s))\myExp{-m(y_A(s))}\right)\differential{s}\\ 
        &{}\quad +\int_{0}^{t}\int_{\mathbb{N}_0}\left(F\left(y_A(s),z-1\right)-F(y_A(s),z)+1\right)^2r_{4}(y_A(s))z\Gamma(\differential{s}\times \differential{z})\\
        &{}\quad +\int_{0}^{t}\int_{\mathbb{N}_0} \lambda \left(F\left(y_A(s),z+1\right)-F(y_A(s),z)\right)^2\Gamma(\differential{s}\times \differential{z})\\
        &{}\quad +\int_{0}^{t}\int_{\mathbb{N}_0}\mu \left(F\left(y_A(s),z-1\right)-F(y_A(s),z)\right)^2 z\Gamma(\differential{s}\times \differential{z}) \\
        &{} = \int^t_0\left(r_1(y_A(s))+r_3(y_A(s))\myExp{-m(y_A(s))}\right)\differential{s}\\ 
        &{}\quad + \int_{0}^{t} \left[ (u_1(y_A(s))-u_2(y_A(s)) + 1)^2 \pi(1) \right. 
        \\
        &{}\quad \quad \quad 
        \left.  + (1- u_2(y_A(s)))^2 (m(y_A(s))-\pi(1)) \right]r_{4}(y_A(s))\differential{s}\\
        &{}\quad + \lambda \int_{0}^{t} \left[(u_2(y_A(s)) - u_1(y_A(s))   )^2 \pi(0) + u_2(y_A(s))^2 \pi([1, \infty) ) \right]\differential{s} \\
        &{}\quad 
        + \mu \int_{0}^{t} \left[ (u_1(y_A(s))-u_2(y_A(s)))^2 \pi(1) + u_2(y_A(s))^2 (m(y_A(s))-\pi(1)) \right] \differential{s} \\
        &{} = \int_0^t \sigma_F(s)\differential{s} \eqcomma 
    \end{align*}
    where $\sigma_F(\cdot)$ is 
    as given
    in the statement of the theorem. Therefore, by the \ac{MCLT} \cite[pp. 339-340]{Ethier:1986:MPC} (see also \cite{Whitt2007Martingale}), we have 
    \begin{align*}
        \tilde{M}^{(n)} \ConvInDist U, 
    \end{align*}
    as $n\to \infty$, where $U$ is a continuous Gaussian martingale with the quadratic variation $\langle U\rangle(t) = \Sigma_F(t)$. In other words, we can write 
    \begin{align*}
        U(t) = \int_{0}^{t} \sqrt{\sigma_F(s)}\differential{B(s)},
    \end{align*}
    for a standard real-valued  Brownian motion $B$. Note that the process $U$ is in $C([0, T], \setOfReals)$ with probability one. 

\paragraph{Convergence of $\frac{1}{\sqrt{n}}(F(Y^{(n)}(t))-F(Y^{(n)}(0)))$}
 Let us now consider the third term in \Cref{eq:main_repres}. Since $F$ is the solution to the Poisson equation, then
   \begin{align*}
&{}\frac{1}{\sqrt{n}}(F(Y^{(n)}(t))-F(Y^{(n)}(0)))\\
&{} =\frac{1}{\sqrt{n}}
 \left(u_1(Y^{(n)}_A(t))\indicator{\{0\}}{Y^{(n)}_B(t)}+u_2(Y^{(n)}_A(t))Y^{(n)}_B(t) \right. \\
 &{}\quad 
 \left. -u_1(Y^{(n)}_A(0))\indicator{\{0\}}{Y^{(n)}_B(0)}-u_2(Y^{(n)}_A(0))Y^{(n)}_B(0)\right)\\
 &{}=\frac{1}{\sqrt{n}}\left(u_1(Y^{(n)}_A(t))\indicator{\{0\}}{Y^{(n)}_B(t)}+u_2(Y^{(n)}_A(t))Y^{(n)}_B(t)-u_1(0)\right).
                 \end{align*}
 Hence, it suffices to check that $\frac{1}{\sqrt{n}}F(Y^{(n)}(t))\ConvInProb 0,\ \text{as}\ n\to\infty$. Note that, for all $t\in[0,T]$, we have 
 \begin{align*}
  \frac{1}{\sqrt{n}}|F(Y^{(n)}(t))|&{}\leqslant  \frac{1}{\sqrt{n}}|u_1(Y^{(n)}_A(t))|+\frac{1}{\sqrt{n}}|u_2(Y^{(n)}_A(t))|| Y^{(n)}_B(t)|.
    \end{align*}
    Hence, we have 
 \begin{align*}
 & \probOf{ \frac{1}{\sqrt{n}}\|F(Y^{(n)}) \|_{C([0,T], \setOfReals)}>\epsilon}  \\ 
 &{}\leqslant \probOf{  \frac{1}{\sqrt{n}}\|u_1(Y^{(n)}_A)\|_{C([0,T], \setOfReals)}+\frac{1}{\sqrt{n}}\|u_2(Y^{(n)}_A)\|_{C([0,T], \setOfReals)}\|Y^{(n)}_B\|_{C([0,T], \setOfReals)}>\epsilon}\eqstop 
    \end{align*}      
  Since 
 $Y^{(n)}_A\ConvInProb y_A$ (as a sequence of stochastic processes in the supremum norm) as $n\to\infty$ by \Cref{thm:relative_compactness} and $u_1,\ u_2$ are continuous functions,
 then due to the continuous mapping theorem
 $\frac{1}{\sqrt{n}}u_1(Y^{(n)}_A)\ConvInProb 0$ and $u_2(Y^{(n)}_A)\ConvInProb u_2(y_A)$ as $n\to\infty.$ We will verify that
 $\frac{1}{\sqrt{n}}\|Y^{(n)}_B\|_{C([0,T], \setOfReals)}\ConvInProb 0$ as $n\to\infty.$ 
To this end, note that for $\beta>0$ and $k\ge 1$, we have 
\begin{align*}
\probOf{\sup_{t\leqslant T}\frac{Y^{(n)}_B(t)}{\sqrt{n}}>k}&{}=\probOf{\tau^{(n)}_{\sqrt{n}k}<T}
=\probOf{e^{-\beta\tau^{(n)}_{\sqrt{n}k}}>e^{-\beta T}}\\
&{}\leqslant  e^{\beta T}\Eof{\myExp{-\beta\tau^{(n)}_{\sqrt{n}k}}}
\leqslant e^{\beta T}\bar{\phi}(n,\sqrt{n}k)\eqcomma 
\end{align*}
by virtue of \Cref{lem:MGF_hitting_times}. Then, in view of \Cref{remark:MGF_hitting_times},  for each $k\geqslant 1$, we have 
$$
\lim_{n\to\infty}\probOf{\sup_{t\leqslant T}\frac{Y^{(n)}_B(t)}{\sqrt{n}}>k} =0.
$$
Therefore, 
$\frac{1}{\sqrt{n}}\absolute{F(Y^{(n)})}\ConvInProb 0\ \text{as}\ n\to\infty$. 

\paragraph{Convergence of the term $\frac{1}{\sqrt{n}}\int^t_0\delta^{(n)}_F(Y^{(n)}(s))\differential{s}$}
It follows from the definition of $\delta_f^{(n)}$ in  \Cref{eq:delta_f_n} that
 \begin{align*}
 &\frac{1}{\sqrt{n}}\int^t_0\delta^{(n)}_F(Y^{(n)}(s))\differential{s}\\
 &{}=\frac{1}{\sqrt{n}}\int^t_0nr_1(Y^{(n)}_A(s))\left[\left(u_1(Y^{(n)}_A(s)+\frac{1}{n})-u_1(Y^{(n)}_A(s))\right)\indicator{\{0\}}{Y^{(n)}_B(s)} \right.\\
&{} \quad \quad \quad 
 \left. +\left(u_2(Y^{(n)}_A(s)+\frac{1}{n})-u_2(Y^{(n)}_A(s))\right)Y^{(n)}_B(s)\right] \differential{s}\\
&{}\quad 
+ \frac{1}{\sqrt{n}}\int^t_0 nr_3(Y^{(n)}_A(s))\indicator{\{0\}}{Y^{(n)}_B(s)}\left[\left(u_1(Y^{(n)}_A(s)-\frac{1}{n})-u_1(Y^{(n)}_A(s))\right)\indicator{\{0\}}{Y^{(n)}_B(s)} \right.\\
&{}\quad \quad \quad 
\left. +\left(u_2(Y^{(n)}_A(s)-\frac{1}{n})-u_2(Y^{(n)}_A(s))\right)Y^{(n)}_B(s)\right] \differential{s}\\
&{}\quad 
+ \frac{1}{\sqrt{n}}\int^t_0 nr_{4}(Y^{(n)}_A(s))Y^{(n)}_B(s)\left[\left(u_1(Y^{(n)}_A(s)-\frac{1}{n})-u_1(Y^{(n)}_A(s))\right)\indicator{\{0\}}{Y^{(n)}_B(s)-1} \right.\\
&{}\quad\quad\quad 
\left. +\left(u_2(Y^{(n)}_A(s)-\frac{1}{n})-u_2(Y^{(n)}_A(s))\right)(Y^{(n)}_B(s)-1)\right]\differential{s} \eqstop \\
 \end{align*}
 Therefore, we have
 \begin{align*}
 & \frac{1}{\sqrt{n}}\left|\int^t_0\delta^{(n)}_F(Y^{(n)}(s))\differential{s}\right| \\
 &{}\leqslant \int^t_0r_1((Y^{(n)}_A(s))\left(|u^{\prime}_1(Y^{(n)}_A(s))|\frac{1}{\sqrt{n}}+\frac{1}{2\sqrt{n}}\left|u^{\prime}_1(Y^{(n)}_A(s)+\frac{1}{n})-u^{\prime}_1(Y^{(n)}_A(s))\right|\right)\differential{s}\\
&{}\quad +\int^t_0r_1((Y^{(n)}_A(s))\\
&{}\quad \quad \quad 
\times\left(|u^{\prime}_2(Y^{(n)}_A)(s)|\frac{1}{\sqrt{n}}+\frac{1}{2\sqrt{n}}\left|u^{\prime}_2(Y^{(n)}_A(s)+\frac{1}{n})-u^{\prime}_2(Y^{(n)}_A)(s))\right|\right)Y^{(n)}_B(s)\differential{s}\\
&{}\quad 
+\int^t_0r_3((Y^{(n)}_A(s))\left(|u^{\prime}_1(Y^{(n)}_A(s))|\frac{1}{\sqrt{n}}+\frac{1}{2\sqrt{n}}\left|u^{\prime}_1(Y^{(n)}_A(s)+\frac{1}{n})-u^{\prime}_1(Y^{(n)}_A(s))\right|\right)\differential{s}\\
&{}\quad 
+\int^t_0r_3((Y^{(n)}_A(s))\\
&{}\quad \quad \quad 
\times\left(|u^{\prime}_2(Y^{(n)}_A)(s)|\frac{1}{\sqrt{n}}+\frac{1}{2\sqrt{n}}\left|u^{\prime}_2(Y^{(n)}_A(s)+\frac{1}{n})-u^{\prime}_2(Y^{(n)}_A)(s))\right|\right)Y^{(n)}_B(s)\differential{s}\\
&{}\quad + \int^t_0r_{4}((Y^{(n)}_A(s))\left(|u^{\prime}_1(Y^{(n)}_A(s))|\frac{1}{\sqrt{n}}+\frac{1}{2\sqrt{n}}\left|u^{\prime}_1(Y^{(n)}_A(s)+\frac{1}{n})-u^{\prime}_1(Y^{(n)}_A(s))\right|\right)\differential{s}\\
&{}\quad 
+\int^t_0r_{4}((Y^{(n)}_A(s))\\
&{}\quad \quad \quad 
\times\left(|u^{\prime}_2(Y^{(n)}_A)(s)|\frac{1}{\sqrt{n}}+\frac{1}{2\sqrt{n}}\left|u^{\prime}_2(Y^{(n)}_A(s)+\frac{1}{n})-u^{\prime}_2(Y^{(n)}_A)(s))\right|\right)|Y^{(n)}_B(s)-1|\differential{s}.
 \end{align*}
Since $\nY_A \ConvInProb y_A$ as $n\to\infty$ by \Cref{thm:relative_compactness}, $\frac{1}{\sqrt{n}}\nY_B \ConvInProb 0$ (both as sequences of stochastic processes in the supremum norm) as $n\to\infty$ from the previous step to show convergence of $\frac{1}{\sqrt{n}}F(Y^{(n)})$, it implies by virtue of the continuous mapping theorem  
$$
\frac{1}{\sqrt{n}}\sup_{t\leqslant T} \left|\int^t_0\delta^{(n)}_F(Y^{(n)}(s))\differential{s}\right| \ConvInProb 0 \text{ as } n\to\infty.
$$
\paragraph{Tightness of the sequence $W_n$}
Denote by
$$
A_n(t)=-\frac{1}{\sqrt{n}}(F(Y^{(n)}(t))-F(Y^{(n)}(0)))
\text{ and }
B_n(t)=\frac{1}{\sqrt{n}}\int^t_0\delta^{(n)}_F(Y^{(n)}(s))\differential{s}.
$$
Then, from \Cref{eq:main_repres}, we have
\begin{align*} 
 |W_n(t)|&{}\leqslant\int^t_0\left(|G^{\prime}(y_A(s))|+\frac{1}{2}|G^{\prime}(Y^{(n)}_A(s))-G^{\prime}(y_A(s))|\right)|W_n(s)|\differential{s} \\
 &{}\quad \quad +\tilde{M}^{(n)}(t)+A(t)+B(t).
\end{align*} 
Applying the Gr\"onwall--Bellman inequality, one can conclude that
\begin{align*}
    \sup_{t\leqslant T}|W_n(t)|& \leqslant\left(\sup_{t\leqslant T}|\tilde{M}^{(n)}(t)|+\sup_{t\leqslant T}|A_n(t)|+\sup_{t\leqslant T}|B_n(t)|\right) \\
     &{}\quad \quad  \times \myExp{\int^T_0\left(|G^{\prime}(y_A(s))|+\frac{1}{2}|G^{\prime}(Y^{(n)}_A(s))-G^{\prime}(y_A(s))|\right)\differential{s}}.
\end{align*}
Since $\predictableVariation{\tilde{M}^{(n)} }(t) \ConvInProb \Sigma_F(t)$, $\sup_{t\leqslant T}|A_n(t)| \ConvInProb 0$,  $\sup_{t\leqslant T}|B_n(t)| \ConvInProb 0$, and $\nY_A \ConvInProb y_A$ as $n\to\infty$, then
by an application of the \ac{BDG} inequality on the martingale, one can verify the compact containment condition (see \Cref{def:compact_containment}) for the sequence $W_n$, \ie, for any $\varepsilon>0$, there exists a positive, finite  constant $k_\varepsilon$ such that
\begin{align}
    \inf_{n}\probOf{\sup_{t\leqslant T}|W_n(t)|\le k_\varepsilon} \ge 1- \varepsilon.
    \label{eq:compact_containment_Wn}
\end{align}

Now, define the module of continuity in $C([0,T], \setOfReals)$ as follows
$$
                    \omega(x,\delta)=\sup_{|t-s|\leqslant\delta}|x(t)-x(s)|\eqcomma 
                    $$
for $x\in C([0,T], \setOfReals)$ and $\delta>0$. 
Then, 
                    \begin{align*}
                    \omega(W_n,\delta)&{}\leqslant\sup_{|t-s|\leqslant\delta} \int^{t\vee s}_{t\wedge s}(|G^{\prime}(y_A(u))|+\frac{1}{2}|G^{\prime}(Y^{(n)}_A(u))-G^{\prime}(y_A(u))|)|W_n(u)|\differential{u}\\  
&{}\quad +\omega(\tilde{M}^{(n)},\delta)+\omega(A_n,\delta)+\omega(B_n,\delta)\\
&{}\leqslant C\delta\sup_{t\leqslant T}|W_n(t)|+\omega(\tilde{M}^{(n)},\delta)+\omega(A_n,\delta)+\omega(B_n,\delta).
 \end{align*} 
 Since
 $\omega(\tilde{M}^{(n)},\delta),\ \omega(A_n,\delta),\ \omega(B_n,\delta)$ converge to zero in probability and in view of \Cref{eq:compact_containment_Wn}, we have that 
 for any $\epsilon>0$
 \begin{align}
 \lim_{\delta\to0}\limsup_{n\to\infty}\probOf{ \omega(W_n,\delta)>\epsilon} =0.
 \label{eq:modulus_continuity}
\end{align}
By virtue of \Cref{eq:compact_containment_Wn} and \Cref{eq:modulus_continuity}, the sequence $W_n$ is tight and hence, relatively compact by Prokhorov's theorem \cite{Billingsley1999Convergence}. 

\paragraph{Identification of the limit}
It is clear that the limit points of the sequence $W_n$  lie in $C([0,T], \setOfReals)$ almost surely.
Choosing a convergent subsequence and passing to the limit in \eqref{eq:main_repres} as $n\to\infty$, we see that a limit point $W$ of the sequence $W_n$ satisfies the stochastic equation
$$
W(t)=\int^t_0 \partial G(y_A(s))W(s)\differential{s} + \int_{0}^{t} \sqrt{\sigma_F(s)}\differential{B(s)}\eqcomma 
$$
or written in the differential form, 
\begin{align*}
    \differential{W(t)} = \partial G(y_A(t))W(t)\differential{t} + \sqrt{\sigma_F(t)}\differential{B(t)}\eqstop 
\end{align*}
Since for each $t>0$, we have 
\begin{align*}
    \int_{0}^{t}\myExp{2 \int_{s}^{t} \partial G(y_A(u))\differential{u} }\sigma_F(s)\differential{s} < \infty \eqcomma
\end{align*}
the above \ac{SDE} admits a unique solution under the assumptions of the theorem. Therefore,  we conclude that the whole sequence $W_n$ converges to the unique solution of the \ac{SDE} in \eqref{eq:clt_sde}. This concludes the proof of the theorem.
\end{proof}

    \section{Discussion}
    \label{sec:discussion}

\subsection*{Practical notes on the model parameterization  
}
    
The formulation of the departure rates of the message queue (queue A) allows modelling a variety of systems that implement different utilization protocols of the buffered entanglement.
The linear dependence of $r_4^{(n)}$ on the number of entanglements originates from the utilization of the buffered entanglements immediately once generated to assist the communication, where each entanglement represents an additional server that goes off once used.
Further, $r_{4}(y_1)$ offers a general model of the rate of message departures per entanglement depending on the utilization protocol of the queued entanglements and the encoding scheme used.
For example, we can apply purification to avoid losing the message at the expense of the number of available entanglements~\cite{elsayed2024fidelity} (and hence the message departure rate per entanglement).
Purification is an essential process in quantum communication that fuses more than one entanglement\footnote{Fusing entanglements is carried out by means of quantum local operations and measurements, which is however not relevant for this work.} to generate a higher fidelity one~\cite{deutsch_quantum_pur} that encodes more information~\cite{bennett1999entanglement} and is more reliable for communication.
Hence, $r_{4}(y_1)$ allows modelling different applications of purification depending on the number of available messages.
We provide in~\Cref{fig: speedup} an example of the assisted service rate $r_4(y)=Mr_3(y_1)y_2$ that comprises a trade-off for the purification protocol. 
The value of $M$ reflects the speedup depending on the average fidelity of entanglements achieved through purification, where a more intense purification, \ie, fusing more entanglements, improves $M$ at the expense of having fewer entanglements $y_2$.
The value of 
$M$ can be determined using existing models~\cite{elsayed2023fidelity,elsayed2024fidelity,davies2023entanglement}.

The message arrival rate $r_1(y_1)$ allows modelling regulated arrivals, while departure rate without entanglement assistance $r_3(y_1)$ allows modelling different classical service models.
The entanglement inter-generation times according to most generation protocols is modelled using a geometric distribution as it requires a constant time to attempt a generation that may fail with some probability $p$~\cite{davies2023tools,elsayed2023fidelity,Dai_Qu_Queuing_delay,Riedel_Nitrogen_vacancy}.
Hence, our model is a continuous time approximation~\cite{ghaderibaneh2023generation,Towsley_stochastic_qu_switch}, especially, since the timescale of generation attempts is very small compared to the system timescale~\cite{pompili2021realization}. 
Although the fidelity decay function is deterministic, modelling the time until discarding an entanglement using an exponential distribution is a common approximation in the literature~\cite{Twosley_ideal_Qu_switch}.
This relies on the interpretation that the randomness of the lifetime of entanglements originates from the randomness of the initial fidelity upon entanglement generation.


\subsection*{Notes on performance-based quantum communication system design}
Using the proposed model we can control the performance of the quantum  communication system by tuning the propensity functions.
Here, we do not particularly tune the classical communication parameters, \ie, regulating the message arrivals $r_1(y_1)$ and the classical service $r_3(y_1)$ as we are mainly interested in the effect of the entanglement assistance. 
The entanglement assistance effect on the performance depends on the entanglement generation rate $\lambda_B$ and the entanglement assisted departure rate $r_4(y_1)$ obtained from the selected entanglement utilization protocols including purification and entanglement generation protocols.

Another interesting parameter to design is the effective service rate, \ie, the steady state expected service rate, which is essentially proportional to the entanglement generation rate and the service rates and is given by
\begin{equation*}
\mu_{\text{eff}}=r_4(y_s) m(y_s) + r_3(y_s) \myExp{-m(y_s)}\eqcomma
\end{equation*}
where $y_s$ is the steady state size of the message queue and $\myExp{-m(y_s)}$ is the stationary probability that the entanglement queue is empty.
For example, consider the case of a classical service rate $r_3(y_1)=\mu_A y_1$, a speedup with a factor M when assisted with entanglements such that  $r_4(y_1)=M\mu_A y_1$ and the arrivals are given by a regulated function $r_1(y_1)=\lambda_A /(1+ y_1)$, where $y_1$ is the size of the message queue.
We obtain the steady state queue size by solving 
\begin{equation*}
     r_1(y_s) -r_3(y_s) \myExp{-m(y_s)} -r_4(y_s)m(y_s)=0\eqcomma
\end{equation*}
with $m(y_s)={\lambda_B}/({M\mu_A y_s+\mu_B})$.
Now, we achieve the desired performance in terms of $y_s$ by tuning the $\lambda_B$ and/or the speedup $M$ by selecting the appropriate entanglement utilization protocol.


\subsection*{Model reparametrization for quantum switch-based graph state distribution}
Our queueing model extends directly to modelling other quantum communication network applications besides the entanglement-assisted communication. 
A major difference is that the message queue in such distributed quantum applications holds requests for consuming entanglements~\cite{elsayed2024fidelity}.
Hence, we set the behaviour of the service offered to requests to be only enabled in the presence of entanglements by setting the propensity function $r_3=0$, while the other propensity functions retain their formulation with $r_4$ parameterizing the service of requests.
For example, this directly applies to applications storing bipartite entanglements between two nodes to use for teleportation, switching or forming graph states (multi-partite entanglements)~\cite{Towsley_stochastic_qu_switch,Dai_Qu_Queuing_delay,dahlberg_qu_Linllayer_protocol}.

    \appendix 
    \section{Additional mathematical details}
    \label{sec:math_details}
    For the sake of completeness, we record here a few basic facts and results that we use in the proofs of the two main theorems, namely, \Cref{thm:relative_compactness} and \Cref{thm:clt}.

    \begin{myDefinition}
        Let $E$ be a complete separable metric space.
        A sequence of $E$-valued stochastic processes $\{U_n\}_{n\ge 1}$ is said to satisfy the \emph{compact containment condition} \citep{Ethier:1986:MPC} if for each $\varepsilon>0$ and $T>0$, there exists a compact set $K_\varepsilon \subseteq E$ such that
        \begin{align*}
            \inf_{n} \probOf{ U_n(t) \in K_\varepsilon, \forall t\in [0, T]} \ge 1-\varepsilon \eqstop
        \end{align*}
        \label{def:compact_containment}
    \end{myDefinition}

    \subsection{\aclp{CTMC}}
    Let $E$ be a countable set. Let $(U(t))_{t\ge 0}$ be a \ac{CTMC}  with paths in $D(\setOfPositiveReals, E)$ with generator $\mathcal{Q}$, formally defined as 
    \begin{align}
        \mathcal{Q}f(x) \defeq \lim_{t\to 0} \frac{\E_x{f(U(t))} - f(x) }{t} \eqcomma 
    \end{align}
    where $\E_x$ denotes the expectation with respect to the probability measure  $\prob_x$ such that $ \prob_x\left({U(0) = x}\right) =1 $, the function $f: E \to \setOfReals$ is bounded, and the limit exists. The proof of the following lemma follows by adapting \cite[Lemma 20.12, p. 33]{Rogers2000DiffusionsVol2}. 

    \begin{myLemma}
        Let $h: \setOfPositiveReals \times E \mapsto \setOfReals$ be a function such that for each $x \in E$, the function $t \mapsto \frac{\partial }{\partial t }h(t, x)$ is continuous. Then, the stochastic process 
        \begin{align*}
            h(t, U(t)) - h(0, U(0)) - \int_{0}^{t} \left(\frac{\partial }{\partial s }h(s, U(s)) + \mathcal{Q}h(s, \cdot)( U(s))  \right)  \differential{s}
        \end{align*}
        is a local martingale. 
        \label{lem:ctmc_Dynkin_timevarying_martingale}
    \end{myLemma}

    \begin{myDefinition}
        A function $h : \setOfPositiveReals \times E \mapsto \setOfReals$ is said to be space-time harmonic for the generator $\mathcal{Q}$ if 
        \begin{align*}
            \frac{\partial }{\partial s }h(s, x) + \mathcal{Q}h(s, \cdot)( x)  = 0 \eqcomma
        \end{align*}
        for all $(s, x) \in \setOfPositiveReals \times E$.
        \label{def:space_time_harmonic}
    \end{myDefinition}

    \begin{myDefinition}
        A function $h : \setOfPositiveReals \times E \mapsto \setOfReals$ is said to be space-time subharmonic (or superharmonic) for the generator $\mathcal{Q}$ if 
        \begin{align*}
            \frac{\partial }{\partial s }h(s, x) + \mathcal{Q}h(s, \cdot)( x)  \ge 0 \quad (\text{or } \le 0)\eqcomma
        \end{align*}
        for all $(s, x) \in \setOfPositiveReals \times E$.
        \label{def:space_time_subsuperharmonic}
    \end{myDefinition}
    
    The following lemma follows from \Cref{lem:ctmc_Dynkin_timevarying_martingale}. See also \cite[Corollary B.5, p. 364]{Robert2003StochasticNetworks}.

    \begin{myLemma}
        Let $h: \setOfPositiveReals \times E \mapsto \setOfReals$ be a function such that for each $x \in E$, the function $t \mapsto \frac{\partial }{\partial t }h(t, x)$ is continuous. If, in addition, the function $h$ is space-time harmonic (or sub- or superharmonic) for the generator $\mathcal{Q}$, then the stochastic process $h(t, U(t))$ is a local martingale (or sub- or supermartingale). 
        \label{lemma:harmonic_martingale}
    \end{myLemma}

    \subsection{Stochastic averaging principle}

    Let $E_1$, and $E_2$ be two complete separable metric spaces, and set $E = E_1\times E_2$. Suppose $\{(U_n, V_n)\}_{n \ge 1}$ be a sequence of stochastic processes with trajectories in $D([0, \infty), E)$, adapted to a filtration $\{\history{t}^{(n)}\}_{t\ge 0}$. Let $\Lambda_n$, an $\mathcal{M}_1([0, \infty), E_2)$-valued random variable, be the occupation measure of the stochastic process $V_n$  defined as 
    \begin{align*}
        \Psi_n (A) \defeq \int_{0}^{t} \indicator{A}{V_n(s)}\differential{s} \eqcomma
    \end{align*}
    for $A \in \borel{E_2}$, and for each $t\ge 0$. 

    The following result, which is the primary instrument used to prove the \ac{FLLN} in \cref{thm:relative_compactness}, is taken verbatim from \cite[Theorem~2.1]{Kurtz1992Averaging}.

    \begin{myTheorem}
    Assume that 
    \begin{enumerate}
        \item The sequence of stochastic processes $\{(U_n)\}_{n \ge 1}$ satisfies the compact containment condition in \Cref{def:compact_containment}. 
        \item The collection $\{V_n(t) \mid t\ge 0, n=1, 2, \ldots\}$ is relatively compact as a collection of $E_2$-valued random variables.
        \item There exists an operator $\mathcal{K} : \effectiveDomain{\mathcal{K}} \subset C_b(E_1, \setOfReals) \to C_b(E_1\times E_2, \setOfReals)$  such that for each $f\in \effectiveDomain{\mathcal{K}}$, there is a stochastic process $\varepsilon_f^{(n)}$ such that the stochastic process 
        \begin{align*}
            f(U_n(t)) - f(U_n(0)) - \int_{0}^{t} \mathcal{K}f(U_n(s), V_n(s)) \differential{s} + \varepsilon_f^{(n)}(t)
        \end{align*}
        is an $\history{t}^{(n)}$-martingale.
        \item The domain $\effectiveDomain{\mathcal{K}}$ of $\mathcal{K}$ is dense in $C_b(E_1, \setOfReals)$ in the topology of uniform convergence on compact sets.
        \item For each $f\in \effectiveDomain{\mathcal{K}}$, and for each $t > 0$, there exists $p>1$ for which
        \begin{align*}
            \sup_{n}\Eof{ \int_{0}^{t} \absolute{\mathcal{K}f(U_n(s), V_n(s))}^p \differential{s}  } < \infty \eqstop 
        \end{align*}
        \item For each $f\in \effectiveDomain{\mathcal{K}}$, and for each $t > 0$, 
        $
            \lim_{n\to \infty} \Eof{ \sup_{s \le t} \absolute{ \varepsilon_f^{(n)}(s) } } = 0 \eqstop
            $
    \end{enumerate}
    Then, the sequence $\{(U_n, \Lambda_n)\}_{n\ge 1}$ is relatively compact in $D([0, \infty), E_1)\times \mathcal{M}_1([0, \infty), E_2)$, and for any limit point $(U, \Lambda)$ of $(U_n, \Lambda_n)$ there exists a filtration $\mathcal{G}_t $ such that
    \begin{align*}
        f(U(t)) - f(U(0)) - \int_{0}^{t} \mathcal{K}f(U(s), v) \Lambda(\differential{s} \times \differential{v})
    \end{align*}
    is a $\mathcal{G}_t$-martingale for each $f\in \effectiveDomain{\mathcal{K}}$.
    \label{thm:stochastic_averaging}
\end{myTheorem}

    
\begin{acronym}[OWL-QN]
	\acro{ABM}{Agent-based Model}
	\acro{BA}{Barab\'asi-Albert}
	\acro{BD}{Birth-Death}
	\acro{BDG}{Burkholder--Davis--Gundy}
	\acro{CDC}{Centers for Disease Control and Prevention}
	\acro{CDF}{Cumulative Distribution Function}
	\acro{CLT}{Central Limit Theorem}
	\acro{CM}{Configuration Model}
	\acro{CME}{Chemical Master Equation}
	\acro{CRM}{Conditional Random Measure}
	\acro{CRN}{Chemical Reaction Network}
	\acro{CTBN}{Continuous Time Bayesian Network}
	\acro{CTMC}{Continuous Time Markov Chain}
	\acro{DSA}{Dynamic Survival Analysis}
	\acro{DTMC}{Discrete Time Markov Chain}
	\acro{DRC}{Democratic Republic of Congo}
	\acro{ER}{Erd\"{o}s-R\'{e}nyi}
	\acro{ESI}{Enzyme-Substrate-Inhibitor}
	\acro{FCLT}{Functional Central Limit Theorem}
	\acro{FLLN}{Functional Law of Large Numbers}
	\acrodefplural{FLLN}[FLLNs]{Functional Laws of Large Numbers}
	\acro{FPT}{First Passage Time}
	\acro{GP}{Gaussian Process}
	\acrodefplural{GP}[GPs]{Gaussian Processes}
	\acro{HJB}{Hamilton–Jacobi–Bellman}
	\acro{iid}{independent and identically distributed}
	\acro{IPS}{Interacting Particle System}
	\acro{KL}{Kullback-Leibler}
	\acro{LDP}{Large Deviations Principle}
	\acro{LLN}{Law of Large Numbers}
	\acrodefplural{LLN}[LLNs]{Laws of Large Numbers}
	\acro{LNA}{Linear Noise Approximation}
	\acro{MAPK}{Mitogen-activated Protein Kinase}
	\acro{MCMC}{Markov Chain Monte Carlo}
	\acro{MCLT}{Martingale Central Limit Theorem}
	\acro{MCT}{Monotone Convergence Theorem}
	\acro{MFPT}{Mean First Passage Time}
	\acro{MGF}{Moment Generating Function}
	\acro{MLE}{Maximum Likelihood Estimate}
	\acro{MM}{Michaelis-Menten}
	\acro{MPI}{Message Passing Interface}
	\acro{MSE}{Mean Squared Error}
	\acro{ODE}{Ordinary Differential Equation}
	\acro{PDE}{Partial Differential Equation}
	\acro{PDF}{Probability Density Function}
    \acro{PDP}{Piecewise Deterministic Process}
	\acro{PGF}{Probability Generating Function}
	\acro{PMF}{Probability Mass Function}
	\acro{PPM}{Poisson Point Measure}
	\acro{PRM}{Poisson Random Measure}
	\acro{psd}{positive semi-definite}
	\acro{PT}{Poisson-type}
	\acro{QSSA}{Quasi-Steady State Approximation}
	\acro{rQSSA}{reversible QSSA}
	\acro{SD}{Standard Deviation}
        \acro{SDE}{Stochastic Differential Equation}
	\acro{SEIR}{Susceptible-Exposed-Infected-Recovered}
	\acro{SI}{Susceptible-Infected}
	\acro{SIR}{Susceptible-Infected-Recovered}
	\acro{SIS}{Susceptible-Infected-Susceptible}
	\acro{sQSSA}{standard QSSA}
	\acro{tQSSA}{total QSSA}
	\acro{TK}{Togashi--Kaneko}
	\acro{WS}{Watts-Strogatz}
	\acro{whp}{with high probability}
\end{acronym}





\section*{Funding}
Olga Izyumtseva was  supported by the British Academy through grant number RaR{\textbackslash}100741, and in part by British Academy, Cara, Leverhulme Trust through grant LTRSF24{\textbackslash}100014. 
The work of W. R. KhudaBukhsh was supported by the Engineering and Physical Sciences Research Council (EPSRC) [grant number EP/Y027795/1].
The work by K. Elsayed was supported by has been funded by the German Research Foundation (DFG) as part of project B4 within the Collaborative
Research Center (CRC) 1053 - MAKI as well as the DFG Cluster of Excellence PhoenixD.

    \bibliographystyle{elsarticle-harv}
    \bibliography{references}

@article{yuan23quantumrouter,
  author  = {Yuan Lee and Eric Bersin and Axel Dahlberg and Stephanie Wehner and Dirk Englund},
  journal = {npj Quantum Inf},
  title   = {A quantum router architecture for high-fidelity entanglement flows in quantum networks},
  year    = {2022},
  volume  = {8},
  number  = {75}
}

@article{fromion2025stochastic,
  title   = {Stochastic Models of Resource Allocation in Chemical Reaction Networks},
  author  = {Fromion, Vincent and Robert, Philippe and Zaherddine, Jana},
  journal = {arXiv preprint arXiv:2511.22252},
  year    = {2025}
}

@article{hashemi2016stochastic,
  title     = {Stochastic averaging and sensitivity analysis for two scale reaction networks},
  author    = {Hashemi, Araz and N{\'u}{\~n}ez, Marcel and Plech{\'a}{\v{c}}, Petr and Vlachos, Dionisios G},
  journal   = {The Journal of chemical physics},
  volume    = {144},
  number    = {7},
  year      = {2016},
  publisher = {AIP Publishing}
}

@inproceedings{bhambay2025optimal,
  title     = {Optimal Scheduling in a Quantum Switch: Capacity and Throughput Optimality},
  author    = {Bhambay, Sanidhay and Vasantam, Thirupathaiah and Walton, Neil},
  booktitle = {Abstracts of the 2025 ACM SIGMETRICS International Conference on Measurement and Modeling of Computer Systems},
  pages     = {139--141},
  year      = {2025}
}

@article{zubeldia2026matching,
  title     = {Matching queues with abandonments in quantum switches: Stability and throughput analysis},
  author    = {Zubeldia, Martin and Jhunjhunwala, Prakirt R and Maguluri, Siva Theja},
  journal   = {Operations Research},
  volume    = {74},
  number    = {1},
  pages     = {339--355},
  year      = {2026},
  publisher = {INFORMS}
}

@article{fustionbasedqcomp,
  author  = {Sara Bartolucci and Patrick Birchall and Hector Bombín and Hugo Cable and Chris Dawson and Mercedes Gimeno-Segovia and Eric Johnston and Konrad Kieling and Naomi Nickerson and Mihir Pant and Fernando Pastawski and Terry Rudolph and Chris Sparrow },
  journal = { Nat Commun},
  title   = {Fusion-based quantum computation},
  year    = {2023},
  volume  = {14},
  number  = {912}
}

@article{azuma2023quantum,
  title     = {Quantum repeaters: From quantum networks to the quantum internet},
  author    = {Azuma, Koji and Economou, Sophia E and Elkouss, David and Hilaire, Paul and Jiang, Liang and Lo, Hoi-Kwong and Tzitrin, Ilan},
  journal   = {Reviews of Modern Physics},
  volume    = {95},
  number    = {4},
  pages     = {045006},
  year      = {2023},
  publisher = {APS}
}

@inproceedings{ghaderibaneh2023generation,
  title        = {Generation and distribution of GHZ states in quantum networks},
  author       = {Ghaderibaneh, Mohammad and Gupta, Himanshu and Ramakrishnan, CR},
  booktitle    = {2023 IEEE International Conference on Quantum Computing and Engineering (QCE)},
  volume       = {1},
  pages        = {1120--1131},
  year         = {2023},
  organization = {IEEE}
}

@inproceedings{elsayed2024fidelity,
  title     = {On the Fidelity Distribution of Purified Link-level Entanglements},
  author    = {Elsayed, Karim S and KhudaBukhsh, Wasiur R and Rizk, Amr},
  booktitle = {IEEE International Conference on Communications (ICC)},
  pages     = {485--490},
  year      = {2024}
}

@article{pompili2021realization,
  title     = {Realization of a multinode quantum network of remote solid-state qubits},
  author    = {Pompili, Matteo and Hermans, Sophie LN and Baier, Simon and others},
  journal   = {Science},
  volume    = {372},
  number    = {6539},
  pages     = {259--264},
  year      = {2021},
  publisher = {American Association for the Advancement of Science}
}

@article{davies2023tools,
  title   = {Tools for the analysis of quantum protocols requiring state generation within a time window},
  author  = {Davies, Bethany and Beauchamp, Thomas and Vardoyan, Gayane and Wehner, Stephanie},
  journal = {arXiv preprint arXiv:2304.12673},
  year    = {2023}
}

@article{davies2023entanglement,
  title   = {Entanglement buffering with two quantum memories},
  author  = {Davies, Bethany and I{\~n}esta, {\'A}lvaro G and Wehner, Stephanie},
  journal = {arXiv preprint arXiv:2311.10052},
  year    = {2023}
}

@article{Riedel_Nitrogen_vacancy,
  title     = {Deterministic Enhancement of Coherent Photon Generation from a Nitrogen-Vacancy Center in Ultrapure Diamond},
  author    = {Riedel, Daniel and S\"ollner, Immo and Shields, Brendan J. and Starosielec, Sebastian and Appel, Patrick and Neu, Elke and Maletinsky, Patrick and Warburton, Richard J.},
  journal   = {Phys. Rev. X},
  volume    = {7},
  issue     = {3},
  pages     = {031040},
  numpages  = {8},
  year      = {2017},
  month     = {Sep},
  publisher = {American Physical Society}
}

@article{Twosley_ideal_Qu_switch,
  title     = {On the exact analysis of an idealized quantum switch},
  author    = {Vardoyan, Gayane and Guha, Saikat and Nain, Philippe and Towsley, Don},
  journal   = {ACM SIGMETRICS Performance Evaluation Review},
  volume    = {48},
  number    = {3},
  pages     = {79--80},
  year      = {2021},
  publisher = {ACM New York, NY, USA}
}

@article{Dai_Qu_Queuing_delay,
  author  = {Dai, Wenhan and Peng, Tianyi and Win, Moe Z.},
  journal = {IEEE Journal on Selected Areas in Communications},
  title   = {Quantum Queuing Delay},
  year    = {2020},
  volume  = {38},
  number  = {3},
  pages   = {605-618},
  doi     = {10.1109/JSAC.2020.2969000}
}

@book{nielsen_QC_QI_book,
  place     = {Cambridge},
  title     = {Quantum Computation and Quantum Information: 10th Anniversary Edition},
  doi       = {10.1017/CBO9780511976667},
  publisher = {Cambridge University Press},
  author    = {Nielsen, Michael A. and Chuang, Isaac L.},
  year      = {2010}
}

@article{Towsley_stochastic_qu_switch,
  title     = {On the stochastic analysis of a quantum entanglement switch},
  author    = {Vardoyan, Gayane and Guha, Saikat and Nain, Philippe and Towsley, Don},
  journal   = {ACM SIGMETRICS Performance Evaluation Review},
  volume    = {47},
  number    = {2},
  pages     = {27--29},
  year      = {2019},
  publisher = {ACM New York, NY, USA}
}

@article{deutsch_quantum_pur,
  title     = {Quantum privacy amplification and the security of quantum cryptography over noisy channels},
  author    = {Deutsch, David and others},
  journal   = {Physical review letters},
  volume    = {77},
  number    = {13},
  pages     = {2818},
  year      = {1996},
  publisher = {APS}
}

@incollection{dahlberg_qu_Linllayer_protocol,
  title     = {A link layer protocol for quantum networks},
  author    = {Dahlberg, Axel and Skrzypczyk, Matthew and Coopmans, Tim and and others},
  booktitle = {Proceedings of the ACM special interest group on data communication},
  pages     = {159--173},
  year      = {2019}
}

@article{cacciapuoti_Caleffi_entang_meets_calssical,
  title     = {When entanglement meets classical communications: Quantum teleportation for the quantum internet},
  author    = {Cacciapuoti, Angela Sara and Caleffi, Marcello and Van Meter, Rodney and Hanzo, Lajos},
  journal   = {IEEE Transactions on Communications},
  volume    = {68},
  number    = {6},
  pages     = {3808--3833},
  year      = {2020},
  publisher = {IEEE}
}

@book{schlosshauer_decoherence,
  title     = {Decoherence: and the quantum-to-classical transition},
  author    = {Schlosshauer, Maximilian A},
  year      = {2007},
  publisher = {Springer Science \& Business Media}
}

@article{elsayed2023fidelity,
  title   = {On the Fidelity Distribution of Link-level Entanglements under Purification},
  author  = {Elsayed, Karim and KhudaBukhsh, Wasiur R and Rizk, Amr},
  journal = {arXiv preprint arXiv:2310.18198},
  year    = {2023}
}

@article{bennett2002entanglement,
  title     = {Entanglement-assisted capacity of a quantum channel and the reverse Shannon theorem},
  author    = {Bennett, Charles H and Shor, Peter W and Smolin, John A and Thapliyal, Ashish V},
  journal   = {IEEE transactions on Information Theory},
  volume    = {48},
  number    = {10},
  pages     = {2637--2655},
  year      = {2002},
  publisher = {IEEE}
}

@article{bennett1999entanglement,
  title     = {Entanglement-assisted classical capacity of noisy quantum channels},
  author    = {Bennett, Charles H and Shor, Peter W and Smolin, John A and Thapliyal, Ashish V},
  journal   = {Physical Review Letters},
  volume    = {83},
  number    = {15},
  pages     = {3081},
  year      = {1999},
  publisher = {APS}
}

@article{hao2021entanglement,
  title     = {Entanglement-assisted communication surpassing the ultimate classical capacity},
  author    = {Hao, Shuhong and Shi, Haowei and Li, Wei and Shapiro, Jeffrey H and Zhuang, Quntao and Zhang, Zheshen},
  journal   = {Physical Review Letters},
  volume    = {126},
  number    = {25},
  pages     = {250501},
  year      = {2021},
  publisher = {APS}
}

@article{bose2000mixed,
  title     = {Mixed state dense coding and its relation to entanglement measures},
  author    = {Bose, Sugato and Plenio, Martin B and Vedral, Vlatko},
  journal   = {Journal of Modern Optics},
  volume    = {47},
  number    = {2-3},
  pages     = {291--310},
  year      = {2000},
  publisher = {Taylor \& Francis}
}

@article{Perry2013Fluid,
  title     = {A Fluid Limit for an Overloaded <i>X</i> Model via a Stochastic Averaging Principle},
  volume    = {38},
  issn      = {1526-5471},
  doi       = {10.1287/moor.1120.0572},
  number    = {2},
  journal   = {Mathematics of Operations Research},
  publisher = {Institute for Operations Research and the Management Sciences (INFORMS)},
  author    = {Perry,  Ohad and Whitt,  Ward},
  year      = {2013},
  month     = {5},
  pages     = {294–349}
}

@article{Laurence2026ScalingMethods,
  title     = {Scaling methods for stochastic chemical reaction networks},
  volume    = {194},
  issn      = {0304-4149},
  doi       = {10.1016/j.spa.2025.104855},
  journal   = {Stochastic Processes and their Applications},
  publisher = {Elsevier BV},
  author    = {Laurence,  Lucie and Robert,  Philippe},
  year      = {2026},
  month     = {4},
  pages     = {104855}
}

@article{Laurence2025AIMD,
  title     = {Stochastic chemical reaction networks with discontinuous limits and AIMD processes},
  volume    = {186},
  issn      = {0304-4149},
  doi       = {10.1016/j.spa.2025.104643},
  journal   = {Stochastic Processes and their Applications},
  publisher = {Elsevier BV},
  author    = {Laurence,  Lucie and Robert,  Philippe},
  year      = {2025},
  month     = {8},
  pages     = {104643}
}

@misc{ganguly2024enzymekinetic,
  doi       = {10.48550/ARXIV.2409.06565},
  author    = {Ganguly,  Arnab and KhudaBukhsh,  Wasiur R.},
  title     = {Statistical inference for a multiscale stochastic model of enzyme kinetics via propagation of chaos},
  publisher = {arXiv},
  year      = {2024},
  copyright = {arXiv.org perpetual,  non-exclusive license}
}

@article{TK2025averaging,
  author  = {Fu, Yi and Kang, Hye-Won and KhudaBukhsh, Wasiur R. and Popovic, Lea and Rempa\l{}a, Grzegorz A. and Williams, Ruth J.},
  title   = {Fragility in a Togashi–Kaneko Stochastic Model with Mutations},
  journal = {SIAM Journal on Life Sciences},
  volume  = {1},
  number  = {2},
  pages   = {202-228},
  year    = {2026},
  doi     = {10.1137/25M1799544}
}

@article{Ganguly2026tQSSA,
  title   = {Asymptotic analysis of the total quasi-steady state approximation for the {M}ichaelis--{M}enten enzyme kinetic reactions},
  journal = {Journal of Mathematical Analysis and Applications},
  volume  = {561},
  number  = {1},
  pages   = {130551},
  year    = {2026},
  issn    = {0022-247X},
  doi     = {https://doi.org/10.1016/j.jmaa.2026.130551},
  author  = {Arnab Ganguly and Wasiur R. KhudaBukhsh}
}

@article{Whitt2007Martingale,
  author     = {Whitt, Ward},
  title      = {Proofs of the martingale {FCLT}},
  journal    = {Probab. Surv.},
  fjournal   = {Probability Surveys},
  volume     = {4},
  year       = {2007},
  pages      = {268--302},
  issn       = {1549-5787},
  mrclass    = {60F17 (60G44)},
  mrnumber   = {2368952},
  mrreviewer = {Dominique\ L\'epingle},
  doi        = {10.1214/07-PS122}
}

@article{Costa2003Poisson,
  title     = {On the Poisson Equation for Piecewise-Deterministic Markov Processes},
  volume    = {42},
  issn      = {1095-7138},
  doi       = {10.1137/s0363012901393523},
  number    = {3},
  journal   = {SIAM Journal on Control and Optimization},
  publisher = {Society for Industrial & Applied Mathematics (SIAM)},
  author    = {Costa,  Oswaldo L. V. and Dufour,  Fran\c{c}ois},
  year      = {2003},
  month     = {1},
  pages     = {985–1001}
}

@article{Glynn1996Liapounov,
  author  = {Glynn, P. W. and Meyn, S. P.},
  journal = {Annals of Probability},
  number  = {2},
  title   = {{A Liapounov bound for solutions of the Poisson equation}},
  volume  = {24},
  year    = {1996},
  pages   = {916-931},
  doi     = {10.1214/aop/1039639370}
}

@book{Robert2003StochasticNetworks,
  author     = {Robert, Philippe},
  title      = {Stochastic networks and queues},
  series     = {Applications of Mathematics (New York)},
  volume     = {52},
  edition    = {French},
  note       = {Stochastic Modelling and Applied Probability},
  publisher  = {Springer-Verlag, Berlin},
  year       = {2003},
  pages      = {xx+398},
  isbn       = {3-540-00657-5},
  mrclass    = {60-02 (60K25 90B22)},
  mrnumber   = {1996883},
  mrreviewer = {Richard\ F.\ Serfozo},
  doi        = {10.1007/978-3-662-13052-0}
}

@book{Rogers2000DiffusionsVol2,
  author    = {Rogers, L. C. G. and Williams, David},
  title     = {Diffusions, {M}arkov processes, and martingales. {V}ol. 2},
  series    = {Cambridge Mathematical Library},
  note      = {It\^o{} calculus,
               Reprint of the second (1994) edition},
  publisher = {Cambridge University Press, Cambridge},
  year      = {2000},
  pages     = {xiv+480},
  isbn      = {0-521-77593-0},
  mrclass   = {60J60 (60G07 60H05 60J25)},
  mrnumber  = {1780932},
  doi       = {10.1017/CBO9781107590120}
}

@book{jacod2003limit,
  title     = {Limit Theorems for Stochastic Processes},
  author    = {Jacod, Jean and Shiryaev, Albert N},
  year      = {2003},
  publisher = {Springer-Verlag Berlin Heidelberg}
}

@article{Crudu2012AAP,
  title     = {Convergence of stochastic gene networks to hybrid piecewise deterministic processes},
  volume    = {22},
  issn      = {1050-5164},
  doi       = {10.1214/11-aap814},
  number    = {5},
  journal   = {The Annals of Applied Probability},
  publisher = {Institute of Mathematical Statistics},
  author    = {Crudu,  A. and Debussche,  A. and Muller,  A. and Radulescu,  O.},
  year      = {2012},
  month     = {10}
}

@book{Applebaum_2009Levy,
  place      = {Cambridge},
  edition    = {2},
  series     = {Cambridge Studies in Advanced Mathematics},
  title      = {L\'evy Processes and Stochastic Calculus},
  publisher  = {Cambridge University Press},
  author     = {Applebaum, David},
  year       = {2009},
  collection = {Cambridge Studies in Advanced Mathematics}
}

@book{IkedaWatanabe2014Stochastic,
  editor    = {Nobuyuki Ikeda and Shinzo Watanabe},
  publisher = {North Holland},
  title     = {{Stochastic Differential Equations and Diffusion Processes}},
  isbn      = {9781483296159},
  year      = {2014}
}

@article{Enger2023Unified,
  title     = {A unified framework for limit results in chemical reaction networks on multiple time-scales},
  volume    = {28},
  issn      = {1083-6489},
  doi       = {10.1214/22-ejp897},
  number    = {none},
  journal   = {Electronic Journal of Probability},
  publisher = {Institute of Mathematical Statistics},
  author    = {Enger,  Timo and Pfaffelhuber,  Peter},
  year      = {2023},
  month     = {1}
}

@article{Hasminskii1966a,
  title     = {On Stochastic Processes Defined by Differential Equations with a Small Parameter},
  volume    = {11},
  issn      = {1095-7219},
  doi       = {10.1137/1111018},
  number    = {2},
  journal   = {Theory of Probability \& Its Applications},
  publisher = {Society for Industrial & Applied Mathematics (SIAM)},
  author    = {Has’minskii,  R. Z.},
  year      = {1966},
  month     = {1},
  pages     = {211–228}
}

@article{Khasminskii1966b,
  title     = {A Limit Theorem for the Solutions of Differential Equations with Random Right-Hand Sides},
  volume    = {11},
  issn      = {1095-7219},
  doi       = {10.1137/1111038},
  number    = {3},
  journal   = {Theory of Probability \& Its Applications},
  publisher = {Society for Industrial & Applied Mathematics (SIAM)},
  author    = {Khas’minskii,  R. Z.},
  year      = {1966},
  month     = {1},
  pages     = {390–406}
}

@article{Kushner2009stochApprox,
  title     = {Stochastic approximation: a survey},
  volume    = {2},
  issn      = {1939-0068},
  doi       = {10.1002/wics.57},
  number    = {1},
  journal   = {WIREs Computational Statistics},
  publisher = {Wiley},
  author    = {Kushner,  Harold},
  year      = {2009},
  month     = dec,
  pages     = {87–96}
}

@book{Kipnis1999Scaling,
  title     = {Scaling Limits of Interacting Particle Systems},
  isbn      = {9783662037522},
  issn      = {0072-7830},
  doi       = {10.1007/978-3-662-03752-2},
  journal   = {Grundlehren der mathematischen Wissenschaften},
  publisher = {Springer Berlin Heidelberg},
  author    = {Kipnis,  Claude and Landim,  Claudio},
  year      = {1999}
}

@book{Billingsley1999Convergence,
  doi       = {10.1002/9780470316962},
  year      = {1999},
  month     = jul,
  publisher = {Wiley},
  author    = {Patrick Billingsley},
  title     = {Convergence of Probability Measures}
}

@inproceedings{Kurtz1992Averaging,
  address   = {Berlin, Heidelberg},
  author    = {Kurtz, Thomas G.},
  booktitle = {Applied Stochastic Analysis},
  editor    = {Karatzas, Ioannis and Ocone, Daniel},
  isbn      = {978-3-540-47017-5},
  pages     = {186--209},
  publisher = {Springer Berlin Heidelberg},
  title     = {Averaging for martingale problems and stochastic approximation},
  year      = {1992}
}

@article{Ball:2006:AAM,
  author    = {Ball, K. and Kurtz, T. G. and Popovic, L. and Rempala, G. A.},
  journal   = {Annals of Applied Probability},
  number    = {4},
  pages     = {1925--1961},
  publisher = {Institute of Mathematical Statistics},
  title     = {Asymptotic analysis of multiscale approximations to reaction networks},
  volume    = {16},
  year      = {2006}
}

@book{Ethier:1986:MPC,
  author    = {Ethier, S. N. and Kurtz, T. G.},
  publisher = {John Wiley \& Wiley},
  title     = {{Markov Processes: Characterization and Convergence}},
  volume    = {282},
  year      = {1986}
}

@article{Kang:2013:STM,
  author    = {Kang, H.-W. and Kurtz, T. G.},
  journal   = {Annals of Applied Probability},
  number    = {2},
  pages     = {529--583},
  publisher = {Institute of Mathematical Statistics},
  title     = {Separation of time-scales and model reduction for stochastic reaction networks},
  volume    = {23},
  year      = {2013}
}

@article{Kang:2014:CLT,
  author    = {Kang, H.-W. and Kurtz, T. G. and Popovic, L.},
  journal   = {Annals of Applied Probability},
  number    = {2},
  pages     = {721--759},
  publisher = {Institute of Mathematical Statistics},
  title     = {Central limit theorems and diffusion approximations for multiscale {M}arkov chain models},
  volume    = {24},
  year      = {2014}
}

\end{document}